\documentclass[12pt]{article}

\usepackage{wrapfig, rotating,color}

\usepackage[utf8]{inputenc} 
\usepackage{textcomp} 
\usepackage{graphicx,epsfig}  
\usepackage{epstopdf} 
\usepackage{flafter}  
\usepackage[numbers]{natbib} 

\usepackage{amsmath,amssymb,amsthm}
\usepackage{bm}  
\usepackage{enumerate}
\usepackage{subfigure}

\newcommand{\IRKA}{\textsf{\small{IRKA}}}

\newcommand{\HinfIRKA}{\textsf{\small{IHA}}}

\newcommand{\MBT}{\textsf{\small{MBT}}}
\newcommand{\BT}{\textsf{\small{BT}}}
\newcommand{\HNA}{{\textsf{\small{HNA}}}}

\newcommand{\Hs }{\ensuremath{H(s)}}
\newcommand{\Hr}{\ensuremath{H_r(s)}}

\newcommand{\boldA}{\ensuremath{\boldsymbol A}}
\newcommand{\boldY}{\ensuremath{\boldsymbol Y}}
\newcommand{\boldX}{\ensuremath{\boldsymbol X}}
\newcommand{\boldE}{\ensuremath{\boldsymbol E}}
\newcommand{\boldZ}{\ensuremath{\boldsymbol Z}}
\newcommand{\boldx}{\ensuremath{\boldsymbol x}}
\newcommand{\boldAr}{\ensuremath{\boldsymbol{A}_r}}
\newcommand{\boldEr}{\ensuremath{\boldsymbol{E}_r}}
\newcommand{\boldb}{\ensuremath{\boldsymbol b}}

\newcommand{\boldbr}{\ensuremath{\boldsymbol {b}_r}}

\newcommand{\boldc}{\ensuremath{\boldsymbol c}}
\newcommand{\boldcr}{\ensuremath{\boldsymbol{c}_r}}

\newcommand{\boldWr}{\ensuremath{\boldsymbol{W}_r}}

\newcommand{\boldP}{\ensuremath{\boldsymbol{P}}}
\newcommand{\boldVr}{\ensuremath{\boldsymbol{V}_r}}
\newcommand{\reals}{\ensuremath{\mathbb{R}}}
\newcommand{\complex}{\ensuremath{\mathbb{C}}}

\newcommand{\IL}{{\mathbb{L}}}
\newcommand{\IM}{{\mathbb{M}}}

\newfont{\Bb}{msbm10 scaled\magstep0}
\def\IR{\mbox {\Bb R}}
\def\IC{\mbox {\Bb C}}
\newcommand{\Hinfty}{\ensuremath{\right\|_ {\mathcal{H}_{\infty}}}}
\newcommand{\Hinfspace}{\ensuremath{\mathcal{H}_{\infty}}}

\newcommand{\boldTheta}{\ensuremath{\boldsymbol{\Theta}}}

\newcommand{\dr}{\ensuremath{d_r}}

\newcommand{\bolde}{\ensuremath{\boldsymbol{e}}}
\newcommand{\boldI}{\ensuremath{\boldsymbol{I}}}

\newcommand{\RHd}{\ensuremath{\{H_r(s,\cdot)\}_{d_r}}}
\newcommand{\SSS}{\textsf{\small SSS}}

\newcommand{\Hinf}{\ensuremath{\mathcal{H}_{\infty}}}
\newcommand{\Htwo}{\ensuremath{\mathcal{H}_{2}}}

\newtheorem{thm}{Theorem}[section]

\theoremstyle{definition}

\theoremstyle{plain}

\title{Interpolatory $\Hinf$  Model Reduction}
\author{Garret Flagg, Christopher Beattie, Serkan Gugercin \\ 
{\small Department of Mathematics, Virginia Tech.} \\
{\small Blacksburg, VA, 24061-0123} \\
{\small \tt e-mail: \{flagg,beattie,gugercin\}@math.vt.edu}
}

\begin{document}
\maketitle
\begin{abstract}
We introduce an interpolation framework for  $\Hinf$ model reduction founded on ideas originating in optimal-$\Htwo$ interpolatory model reduction, realization theory,  and complex Chebyshev approximation.  By employing a Loewner ``data-driven" framework within each optimization cycle,  large-scale $\Hinf$ norm 
calculations can be completely avoided. Thus, we are able to formulate a method that remains effective in large-scale settings with the main cost dominated by sparse linear solves.
Several numerical examples 
 illustrate that our approach will produce high fidelity  reduced models  consistently exhibiting better $\Hinf$ performance than those produced by balanced truncation; these models often are as good as (and occasionally better than) those models produced by optimal Hankel norm approximation.  In all cases, these reduced models are produced at far lower cost than is possible either with balanced truncation or optimal Hankel norm approximation.
\end{abstract}

\section{Introduction}
The need for high accuracy mathematical models in problems involving simulation and control often results in dynamical systems described by a large number of differential equations.  Working with such large-scale systems can easily place overwhelming demands on computational resources, a problem which model reduction methods seek to alleviate by approximating the original model with another model consisting of far fewer (but carefully crafted) differential equations.  Strategies for carrying out this approximation should be both efficient and accurate.  For an overview of model reduction, see \cite{antB}.

We consider here single-input/single-output (SISO) linear dynamical systems given in state-space form as:
\begin{eqnarray}
\label{ltisystemintro}
\boldE \dot{\boldx}(t) = \boldA \boldx(t) + \boldb\,u(t),\qquad
y(t) = \boldc^T\, \boldx(t) + d\, u(t),
\end{eqnarray}
where $\boldE,\, \boldA \in \IR^{n \times n}$, $\boldb,\boldc \in \IR^{n}$ and $d \in \IR$.  $\boldE$ is assumed to be nonsingular throughout, although our approach extends without difficulty to cases where $\boldE$ is singular so long as $\mbox{\textsf{nullity}}(\boldE)=\mbox{\textsf{nullity}}(\boldE^2)$ (that is, provided that $0$ is a nondefective eigenvalue for $\boldE$).  
 $\boldx(t)\in \IR^n$ are the \emph{states}; $u(t)\in \IR$ is the \emph{input}; and $y(t)\in \IR$ is the \emph{output} of the dynamical system 
in (\ref{ltisystemintro}).
The \emph{transfer function} of the system is 
$$
H(s) =\boldc^T(s\boldE-\boldA)^{-1}\boldb+d,
$$
defined for $s\in\IC$.  In accord with standard convention, 
we  denote both the system and its transfer function by $H(s)$.
We assume that $H(s)$ is both controllable and observable. 
The \emph{order} of $H(s)$ is the number of poles it possesses, counting multiplicity.  
Since $\boldE$ is nonsingular, all poles of $H(s)$ are finite and since $H(s)$ is both controllable and observable,  
the order of $H(s)$ is identical to the dimension, $n$, of the state vector $\boldx$ in  (\ref{ltisystemintro}).

 We denote by $\Hinf^k$, the set of rational functions of order at most $k$ which are bounded and analytic in the closed right half plane in $\mathbb{C}$. We assume in all that follows that $H \in \Hinf^n$.
 The $\Hinf$ norm of $H$ is defined as
\begin{equation} \label{hinfDef}
\left\| H \right\|_{\Hinf}=
\max_{\omega\in {\mathbb R}}\mid H(\jmath \omega) \mid.
\end{equation}
  If  the input function, $u(t)$, is square integrable:
 $\displaystyle \|u\|_{L^2}=\left(\int_0^\infty \mid u(t) \mid ^2\,dt \right)^{\frac{1}{2}} <\infty$, then
 the output function, $y(t)$,  of (\ref{ltisystemintro}) will be square integrable as well; 
 $u$ and $y$ have Fourier transforms  $\widehat{u},\,\widehat{y}\in L^2(\IR)$ that are related according to
 $\widehat{y}(\omega)=H(\jmath \omega)\widehat{u}(\omega)$.  One immediately observes
 that  the $\Hinf$ norm defined in (\ref{hinfDef}) is an $L^2$-induced
operator norm of the underlying convolution operator mapping $u\mapsto y$.
 
Our goal is  to  construct another system
\begin{eqnarray} \label{redsysintro}
\boldEr \dot{\boldx}_r (t) = \boldA_r \boldx_r (t) + \boldb_r u(t),\qquad
y_r(t)  =  \boldcr^T \boldx_r (t)  + d_r u(t)
\end{eqnarray}
of much smaller order $r \ll n$, with  $\boldEr,\, \boldA_r \in \IR^{r \times r}$,
$\boldb_r,\,\boldcr \in \IR^{r}$, and $d_r\in \IR$ determined  so
 that $y_r$ approximates $y$ uniformly well over all 
 $u\in  L^2(\IR^+)$, in an appropriate sense.

Toward this end, define a \emph{reduced transfer function} associated with (\ref{redsysintro}) as
$H_r(s) = \boldcr^T(s\boldEr-\boldA_r)^{-1}\boldb_r+d_r$.
 Then, for any $u\in L^2(\IR^+)$, 
\begin{equation} \label{eqn:L2induced}
\| y- y_{r} \|_{L^2} \leq  \| H- H_{r} \|_{\Hinf} \| u \|_{L^2}.
\end{equation}
The \emph{error transfer function}, $H(s)-H_r(s)$, may be rewritten as 
$$
H(s) - H_r(s) = \boldc^T (s\boldE-\boldA)^{-1}\boldb - \left[ \boldcr^T(s\boldEr-\boldA_r)^{-1}\boldb_r+ (d_r-d)\right].
$$
Evidently, a non-zero $d$-term in the original model may be absorbed into a reduced-order model by assigning $d_r\mapsto (d_r-d)$.  
This allows us to assume in all that follows that $d=0$ in the original model without any loss of generality. 

In view of  (\ref{eqn:L2induced}),  the output error magnitude, $\|y- y_{r}\|_{L^2}$, may be made uniformly small over all bounded inputs $u$ (say, with $\|u\|_{L^2}\leq 1$) 
if we find a reduced system, $H_r$, that makes  $ \| H- H_{r} \|_{\Hinf}$ small. 
 This leads naturally to the 
\emph{optimal $\Hinf$ model reduction problem}:  

For a given $H \in \Hinf^n$
and reduction order $r < n$, find $H_r^{\star} \in \Hinf^r$,
that solves 
 \begin{equation} \label{optimalHinf}
 \min \limits_{\hat{H}_r \in \Hinf^r} \left\| H-\hat{H}_r \Hinfty.
\end{equation}

This problem is an active area of research \cite{antoulas2005hna}.   We develop a methodology for approximating solutions to (\ref{optimalHinf})  that remain effective in large-scale settings,
settings where the original state-space dimension, $n$, could be on the order of $10^5$ or more, for example.   Most  methods known to us will be intractable even for modest system order, say, on the order of a few thousand (with the exception of balanced truncation, see the discussion below).

 Kavrano{\u{g}}lu and Bettayeb  \cite{Kavvy} showed that (\ref{optimalHinf}) can be 
converted into an optimal Hankel norm approximation problem for a special imbedded system with augmented input and output mappings.
However, this approach is infeasible in practice since knowledge both of the minimum of (\ref{optimalHinf}) as well the imbedded system is required. 
As noted in \cite{Kavvy}, this information is available (or computationally accessible) only in very special cases. 

Several methods to solve (\ref{optimalHinf}) that utilize \emph{linear matrix inequality} (LMI) frameworks have been presented as well; see, for example, \cite{grigoriadis1995ohm,kavranoglu1993zoh,kavranoglu1993had,helmersson2002mru,varga2001fas} and references therein.  
These approaches rapidly become computationally intractable with increasing state space dimension.  Indeed, published examples illustrating LMI-based methods in \cite{grigoriadis1995ohm,kavranoglu1993zoh,kavranoglu1993had,helmersson2002mru,varga2001fas} all had order less than $n=10$. 

The most common practical methods for obtaining satisfactory  $\Hinf$ reduced models are
Gramian-based methods such as \emph{balanced truncation} (\BT) \cite{moore1981principal,mullis1976synthesis}
and \emph{optimal Hankel norm approximation methods} (\HNA) \cite{Gloverpaper}.  
Both approaches are known to yield small approximation errors in the $\Hinf$ norm \cite{gugercin2004asurvey,antB}, though neither generally is capable of producing globally optimal 
solutions to (\ref{optimalHinf}).   Both approaches also remain computationally feasible for modest state space dimension (perhaps a few thousand), but significantly larger state space dimension still presents challenges.  \HNA~requires an all-pass dilation of the full-order model followed by a full eigenvalue decomposition. These are dense matrix operations, having a complexity growing with $\mathcal{O}(n^3)$, which generally limits problem sizes to a few thousand.  Notably,  \cite{benner2004coh} was able to extend \HNA~to dynamical systems having state space dimension of $\mathcal{O}(10^4)$, using state-of-the-art numerical techniques tailored to high performance computer architectures. 

The situation for \BT~is a bit better.  \BT~has been applied to systems of order $\mathcal{O}(10^5)$ by solving the underlying Lyapunov equations {\it iteratively} using ADI-type algorithms; see, for example, \cite{gugercin2003amodified,stykel2004gbm,penzl200clr,benner2003statespace,sorensen2002thesylvester,heinkenschloss2008btm} and  references therein.  


We describe here a different model reduction methodology that can be  applied effectively even for very large state space dimension, yielding reduced models typically having  smaller $\Hinf$ errors than either \BT~or \HNA\, moreover, at significantly lower cost. 
Towards this end, we  present a new framework for the $\Hinf$ approximation problem using  interpolatory model reduction. By connecting ideas from interpolatory $\Htwo$ model reduction \cite{H2}, realization theory \cite{paramD},  and complex Chebyshev approximation \cite{Trefsufficient}, we  develop an  interpolation-based method  for  $\Hinf$-approximation that remains numerically efficient in large-scale settings.  The main cost of our approach involves the solution of sparse linear systems. We  demonstrate that our approach can yield reduced-order systems having  $\Hinf$ errors  which are often half that of \BT, and very close to (indeed, sometimes better than) that of  \HNA.   For symmetric systems, our method typically produces reduced-order systems with $\Hinf$-errors that are near the theoretical best possible of (\ref{optimalHinf}).

The rest of the paper is organized as follows: we close this section with a brief review of interpolatory model reduction and related approaches for 
solving the optimal $\Htwo$ approximation problem.  We introduce our method in \S \ref{sec:hinf_general}  
and illustrate its effectiveness via several numerical examples in \S \ref{sec:examples}. 

\paragraph{Interpolatory model reduction.} \label{sec:interpolation}
Given a dynamical system \Hs\ and  a set of points $\{s_1,\,s_2,\,\ldots,\,s_k\}\subset \IC$,   \emph{interpolatory model reduction}  produces a dynamical system \Hr\ such that \Hr\ interpolates \Hs\ together with a prescribed number of derivatives at the points $\{s_1,\,s_2,\,\ldots,\,s_k\}$.  Although this  is posed as a rational interpolation problem,  the construction of a solution may be accomplished with a variety of rational Krylov subspace projection techniques.  Rational interpolation via projection was first proposed  by Skelton \emph{et al}. \cite{Skeltonlate,Skelt1,Skelt2}.  Later, Grimme \cite{Grim} showed how to construct a reduced-order (Hermite) interpolant using a method of Ruhe.  
%
%
\begin{thm}[Grimme \cite{Grim}]\label{thm:rationalkrylov}
Given $\Hs=\boldc^T(s\boldE-\boldA)^{-1}\boldb$ and a point-set $\mathcal{S}_1 \subset \complex$ containing $r$ distinct points: $\mathcal{S}_1=\{s_1, \ldots, s_{r}\} $, let
\begin{equation} \label{eqn:VrWr}
\boldsymbol{V}_r=\lbrack (s_1\boldE -\boldsymbol{A})^{-1}\boldsymbol{b} \dots  (s_r\boldE -\boldsymbol{A})^{-1}\boldsymbol{b}\rbrack
\qquad
\boldsymbol{W}_r^T=
\begin{bmatrix}
\boldc^T(s_{1}\boldE-\boldsymbol{A})^{-1}\\
\vdots \\
\boldc^T(s_{r}\boldE-\boldsymbol{A})^{-1}
\end{bmatrix}.
\end{equation}

Define a reduced-order model $H_r^{0}(s) = \boldcr^T(s \boldEr - \boldAr)^{-1}\boldbr$, where $d_r=0$ and
\begin{equation}\label{defineVW}
\boldEr= \boldsymbol{W}_r^T\boldE\boldsymbol{V}_r,~~\boldsymbol{A}_r=\boldsymbol{W}_r^T\boldsymbol{AV}_r,~~ \boldsymbol{b}_r=\boldsymbol{W}_r^T\boldsymbol{b},~~{\rm and}~~
 \boldcr^T= \boldc^T\boldsymbol{V}_r.
 \end{equation}
Then $H(s_i)=H_r^{0}(s_i)$ and $H'(s_k)=H_r^{0\,\prime}(s_k)$, for $ i=1, \dots 2r$  where $\phantom{}'$ denotes the derivative with respect to the frequency parameter, $s$.
\end{thm}
Higher-order derivatives can be matched similarly; for details, see 
\cite{Grim,Ant2010imr}. 

\paragraph{$\Htwo$-optimality conditions.}\label{sec:irka}
Theorem \ref{thm:rationalkrylov} gives explicitly computable conditions that will yield a reduced-order model satisfying $2r$ interpolation conditions; one need only solve $2r$ linear algebraic systems to form the columns of $\boldsymbol{V}_r$ and $\boldsymbol{W}_r$.   Notably, Theorem \ref{thm:rationalkrylov} carries no hint of how best to choose these interpolation points.     A method for determining interpolation points that leads to reduced-order models that are (locally) optimal with respect to the \Htwo\ error was developed in \cite{H2} and will be a point of departure for the approach we propose here.  

For a SISO dynamical system \Hs, the \Htwo\ norm is defined as
\begin{equation}
\left \| H \right \|_{\Htwo} = \left(\frac{1}{2\pi} \int_{-\infty}^\infty \mid H(\jmath \omega) \mid^2 d \omega \right)^{1/2}.
\end{equation}
Then, for a given full-order model $H(s)$,
and selected reduction order $r < n$, 
the \emph{optimal $\Htwo$ approximation problem} seeks a reduced model, ${\displaystyle H_r^{0}(s)}$,
that solves 
 \begin{equation} \label{optimalHtwo}
 \min \limits_{\hat{H}_r^{0} \in \Htwo^r}   \left\| H-\hat{H}_r^{0} \right\|_{\Htwo}.
\end{equation}
The approximation problem (\ref{optimalHtwo}) has been studied extensively; 
see, for example,
\cite{meieriii1967approximation},\cite{wilson1970optimum},\cite{H2}, \cite{spanos1992anewalgorithm}, \cite{vandooren2008hom}, \cite{gugercin2005irk}, \cite{beattie2007kbm}, \cite{beattie2009trm},\cite{zigic1993contragredient} and  references therein.
First-order necessary conditions for \Htwo-optimal approximation may be formulated in terms of interpolation conditions:
\begin{thm}[\cite{meieriii1967approximation,H2}] \label{thm:h2optimal}
Given a full-order model \Hs, let $H_r^{0}(s)$ be an \Htwo-optimal reduced order model of order $r$, with simple poles $\hat{\lambda}_1,\dots,\hat{\lambda}_r$. Then,
\begin{equation} \label{eqn:h2cond}
H(-\hat{\lambda}_i)=H_r^{0}(-\hat{\lambda}_i)~~~{\rm and}~~~
H'(-\hat{\lambda}_i)=H_r^{0\,\prime}(-\hat{\lambda}_i)~~~{\rm for}~~~
 i=1,\dots,r,
\end{equation}
where $\phantom{}'$ denotes differentiation with respect to the frequency parameter, $s$.
\end{thm}

These necessary conditions characterize the \Htwo-optimal reduced order model as a rational Hermite interpolant matching the full-order transfer function and its derivative at mirror images of the reduced-order system poles. 
This can be accomplished with the help of Theorem \ref{thm:rationalkrylov} once the poles of $H_r^{0}(s)$ are known.  Of course,  the poles of $H_r^{0}(s)$ are not known \emph{a priori}, but they can be computed iteratively using the \emph{Iterative Rational Krylov Algorithm} ({\IRKA}) developed by Gugercin \emph{et al.} \cite{H2}.  

\begin{figure}[ht]
\begin{center}
    \framebox[6.25in][t]{
    \begin{minipage}[c]{5.75in}
    {
\textbf{\textsf{\small{Algorithm}}~\IRKA}.  \emph{Iterative Rational Krylov Algorithm \cite{H2}} 

\quad Given a full-order $H(s)$, a reduction order $r$,  and convergence tolerance $\mathsf{tol}$.
\begin{enumerate}
\item Make an initial selection of interpolation points $s_i$, for $i=1,\dots,r$ that is closed under complex conjugation.
\item Construct \boldVr\ and \boldWr\ as in (\ref{eqn:VrWr}) with $s_{i+r}=s_i$ for $i=1,\ldots,\,r$
\item while (relative change in $\{s_i\} < \mathsf{tol} $)
\begin{enumerate}[a)]
\item $\boldEr=\boldWr^T \boldE\boldVr$ and $\boldAr= \boldWr^T \boldA\boldVr$
\item Solve $r\times r$ eigenvalue problem $\boldAr\mathbf{u}=\lambda\boldEr\mathbf{u}$.
\item Assign $s_i= s_{i+r} \leftarrow -\lambda_i(\boldAr,\boldEr)$ for $i=1,\dots ,r$.
\item Update \boldVr\ and \boldWr\  as in (\ref{eqn:VrWr}) using new $\{s_i\}$.
\end{enumerate}
\item $\boldEr=\boldWr^T \boldE\boldVr$, $\boldAr= \boldWr^T \boldA\boldVr$, $\boldbr=\boldWr^T\boldb$, $\boldcr^T=\boldc^T\boldVr$,  and $d_r= 0$.
\end{enumerate}
 }
    \end{minipage}
    }
    \end{center}
  \end{figure}

\IRKA\ is a fixed point iteration that in the SISO case typically  exhibits rapid convergence  to a local minimizer of the \Htwo-optimal model reduction problem. 
Sparsity in $\boldE$ and $\boldA$ can be well-exploited in the linear solves of Steps 2 and 3c and \IRKA\ has been remarkably successful in producing high fidelity reduced-order approximations in large-scale settings;  it has
been  applied successfully in finding $\Htwo$-optimal reduced models for systems of high order (e.g., $n>160,000$, see \cite{KRXC08}).  For details on the algorithm, we refer to the original source \cite{H2}. 

\section{An interpolatory approach for $\Hinf$ approximation}  \label{sec:hinf_general}
 
The principal result that defines the character of our approach was provided by Trefethen in 
\cite{Trefsufficient}; it is an analog to the Chebyshev Equioscillation Theorem.  

\begin{thm}[Trefethen \cite{Trefsufficient}] \label{thm:trefsufficient}
Suppose $H(s)$ is a transfer function associated with a dynamical system as in (\ref{ltisystemintro}).  Let $H_r^{\rm opt}(s)$ be an optimal $\Hinf$ approximation to $H(s)$ (i.e., a solution to (\ref{optimalHinf}) ) and let \Hr\  be any $r^{\rm th}$ order stable approximation to $H(s)$ that interpolates $H(s)$ at $2r+1$ points in the open right half plane. 
Then
$$
\min_{\omega \in \IR} | H(\jmath \omega) - H_r(\jmath \omega)| \leq  \| H - H_r^{\rm opt} \|_{\Hinf} \leq \| H - H_r \|_{\Hinf} 
$$
In particular,  if
$|H(\jmath \omega)-H_r(\jmath \omega)|=\mathsf{const}$ for all $\omega\in\IR$ 
then \Hr\ is itself an optimal $\Hinf$ approximation to $H(s)$.	
\end{thm}

%
%
One  sees from this that a good $\Hinf$ approximation will be obtained when the modulus of the error, $|H(s)-H_r(s)|$, is nearly constant as $s=\jmath \omega$ runs along the imaginary axis. 
We select $2r+1$ interpolation points in the open right half plane that will induce this.  By utilizing Theorem \ref{thm:rationalkrylov}, we may locate $2r$ interpolation points in the right half plane as we like, producing an interpolating reduced-order system, $H_r^{0}(s)$.  Also,  $H(\infty)=H_r^{0}(\infty)=0$ so we can exploit the freedom in choosing  $d_r$ to move the $(2r+1)^{\rm st}$ interpolation point from $\infty$ into the open right half-plane.  Note that the straightforward construction, $H_r(s,d_r)= H_r^{0}(s)+d_r$,  creates a reduced-order model that has \emph{all} interpolation points depending on $d_r$.   We prefer a different formulation that uncouples the $2r$ interpolation points associated with 
$H_r^{0}(s)$ from the influence of $d_r$.    The construction that accomplishes this was introduced in \cite{paramD} (see also \cite{Struct} for a formulation close to what we use here).
\begin{thm} \label{thm:sameprojection}
Given $H(s) = \boldc^T(s \boldI -\boldA)^{-1}\boldb$ and a point-set $\mathcal{S}\ \subset \complex$ containing $r$ distinct points: $\mathcal{S}=\{s_1, \ldots, s_{r}\} $, let 
$\boldsymbol{V}_r$, $\boldsymbol{W}_r$, $\boldEr$, $\boldA_r$, $\boldb_r$, and $\boldcr$ be defined as in
Theorem \ref{thm:rationalkrylov}. For any given $d_r \in \IR$, define a new reduced-order system 
\begin{equation} 
H_r(s,d_r) = (\boldc_r-d_r\bolde)^T(s\boldEr-\boldA_r - d_r\bolde \bolde^T)^{-1}(\boldbr-d_r\bolde) + d_r \label{Htilde}
\end{equation}
where $\bolde$ denotes a vector of ones.  
Define auxiliary reduced systems: 
\begin{align*}
H_r^{0}(s)=\boldcr^T(s \boldEr -\boldA_r)^{-1}\boldb_r, \quad & \quad
G_1(s)=\bolde^T(s \boldEr -\boldA_r)^{-1}\boldb_r, \\
G_2(s)=\boldcr^T(s \boldEr -\boldA_r)^{-1}\bolde, \quad \mbox{and} & \qquad 
G_3(s)=\bolde^T(s \boldEr -\boldA_r)^{-1}\bolde.
\end{align*}
Then 
\begin{equation} \label{drformula}
H_r(s,d_r)=H_r^{0}(s)+d_r \frac{(G_1(s)-1)(G_2(s)-1)}{1-d_r G_3(s)}
\end{equation}
 and for all $d_r \in \IR$
$$
H(s_i)=H_r^{0}(s_i)=H_r(s_i,d_r)~~\mbox{ and }~~H'(s_i)=H_r^{0\,\prime}(s_i)=H_r'(s_i,d_r)\quad \mbox{ for }\quad i=1,\ldots,r
$$
where $\phantom{}'$ denotes the derivative with respect to the frequency parameter, $s$.
\end{thm} 

\noindent\textsc{Proof:}\ The expression (\ref{drformula}) follows from (\ref{Htilde}) with 
straightforward manipulations that begin with the Sherman-Morrison formula: 
$$
(s\boldEr-\boldA_r - d_r\bolde \bolde^T)^{-1}=(s\boldEr-\boldA_r)^{-1}+\frac{d_r}{1-d_r G_3(s)}(s\boldEr-\boldA_r)^{-1}\bolde \bolde^T(s\boldEr-\boldA_r)^{-1}
$$

Define $\boldP(s)=\boldVr(s\boldEr-\boldA_r)^{-1}\boldWr^T(s\boldE-\boldA)$.  Observe that $\boldP(s)$ is a (skew) projection onto $\mbox{\textsf{Ran}}(\boldVr)$ for any $s\in\IC$ for which it is well-defined and so we have
\begin{align*}
\boldVr\bolde_k=&\boldP(s_k)\boldVr\bolde_k=\left[\boldVr(s_k\boldEr-\boldA_r)^{-1}\boldWr^T(s_k\boldE-\boldA)\right](s_k\boldE-\boldA)^{-1}\boldb\\
&=\boldVr(s_k\boldEr-\boldA_r)^{-1}\boldWr^T\boldb=\boldVr(s_k\boldEr-\boldA_r)^{-1}\boldb_r
\end{align*}
Linear independence of the columns of $\boldVr$ then implies $\bolde_k=(s_k\boldEr-\boldA_r)^{-1}\boldb_r$ and thus, $G_1(s_k)=1$, for $k=1,\,2,\,\ldots,\,r$. 
 A similar argument yields $G_2(s_k)=1$, for $k=1,\,2,\,\ldots,\,r$. Taken together, we get that $H_r^{0}(s_i)=H_r(s_i,d_r)$ for $i=1,\,2,\,\ldots,\,r$. Likewise, 
{\small 
 $$
 H_r'(s,d_r)-H_r^{0\, \prime}(s)= \frac{d_r^2 G_3'(s)}{(1-d_r G_3(s))^2}(G_1(s)-1)(G_2(s)-1)
+ d_r \frac{G_1'(s)(G_2(s)-1)}{1-d_r G_3(s)}
 +d_r \frac{(G_1(s)-1)G_2'(s)}{1-d_r G_3(s)}
 $$
}
so that $H_r^{0\,\prime}(s_i)=H_r'(s_i,d_r)$ as well. $\Box$

%
%
%
\vspace{1ex}
\noindent
 Let \RHd\ denote the set of all transfer functions $H_r(s,d_r)$ with $\dr$ ranging over $\reals$. The freedom we have in choosing \dr\  is significant to us for at two reasons.  First, \RHd\ is a parameterization of the set of all proper rational functions of degree $r$ having real coefficients that satisfy the same interpolation constraints as $H_r^{0}(s)$ (see e.g., \cite{paramD}). Second, it is now possible to construct reduced-order models of order $r$ satsfying $2r+1$ interpolation conditions, which is an essential step towards constructing reduced-order models that are optimal in the \Hinfspace-norm.  Since $H_r(s,d_r)$ interpolates $H(s)$ at $s_1,\ldots,s_{2r}$ for any $d_r$, one could select an additional (real) interpolation point, $s_{2r+1}>0$ and 
directly calculate from (\ref{drformula}) the value of $d_r$ that enforces $H_r(s_{2r+1},d_r) = H(s_{2r+1})$:
$$
d_r=\frac{H(s_{2r+1})-H_r^{0}(s_{2r+1})}{(G_1(s_{2r+1})-1)(G_2(s_{2r+1})-1)+G_3(s_{2r+1})(H(s_{2r+1})-H_r^{0}(s_{2r+1}))}.
$$
We avoid the necessity of explicitly selecting $s_{2r+1}$.  Instead, as discussed below, $d_r$ will be chosen directly to decrease the $\Hinf$ error. 

\subsection{An algorithm  for $\Hinf$ approximation}\label{sec:heuristic1}
%

It has been observed (e.g., see \cite{H2},\cite{Ant2010imr}) that $\Htwo$ optimal interpolation points produced by \IRKA~yield reduced  models that are not only  (locally) $\Htwo$-optimal but  frequently also are high-fidelity $\Hinf$ approximations  to the original system.   Indeed, $\Htwo$-optimal models produced by \IRKA\ yield $\Hinf$ error norms that are comparable 
to that of  \BT\  and sometimes are even better. 
 Therefore, our approach begins with \textsf{\small{Algorithm}}~\IRKA\ to obtain $2r$  interpolation points (counting multiplicity -- Hermite interpolation at $r$  distinct points) determining an $\Htwo$-optimal reduced model, $H_r^0(s)$.   This  choice for $H_r^0(s)$ defines a family of approximations parameterized by $d_r$, \RHd. We then proceed by (approximately)
minimizing the $\Hinf$ error with respect variations in $d_r$.   These steps are summarized below: 


\begin{figure}[hhh]
\begin{center}
    \framebox[6.25in][t]{
    \begin{minipage}[c]{5.80in}
    {
\textbf{\textsf{\small{Algorithm}}~\HinfIRKA}.  \emph{Interpolatory $\Hinf$ Approximation:} 

\quad Given a full-order model, $H(s)$, and reduction order, $r$.
\begin{enumerate}
\item Apply \textsf{\small{Algorithm}}~\IRKA~to compute $2r$ $\Htwo$-optimal interpolation points and an associated $\Htwo$-optimal reduced model, $H_r^0(s)$.
\item Find 
$\displaystyle d_r^\star=\arg \min \limits_{d_r \in \reals} \left \| H-H_r \Hinfty$ where $H_r=H_r(s,d_r)$ is defined in (\ref{drformula}).
\item Construct the final $\Hinf$ approximant as $\displaystyle H_r^{\star}(s)=H_r(s,d_r^\star)$. 
\end{enumerate}
       }
    \end{minipage}
    }
    \end{center}
  \end{figure}

The distribution of interpolation points obtained in Step 1 (as an outcome of  $\Htwo$-optimal approximation) yields very effective $\Hinf$ approximants as well.  We therefore wish to preserve these interpolation points using Theorem \ref{thm:sameprojection} while varing the $\dr$ parameter in such a way as to drive down the $\Hinf$ error - centering the error curve about the origin in the process.   ${\displaystyle H_r^\star(s)}$  will  denote an $\Hinf$ approximant having 1) an $\Htwo$ optimal pole distribution (from Step 1 of \textsf{\small{Algorithm}}~\HinfIRKA) and 2) an optimally chosen  $d_r^\star$ (from Step 2).



\subsection{Efficient  Implementation of  Step 2 of \textsf{{Algorithm IHA}}}  \label{loewner}
The major contributions to the cost of \HinfIRKA~come from linear solves arising in Step 1 (from \IRKA) and the $\Hinf$ norm evaluations required in solving the  (scalar) nonlinear optimization problem  in Step $2$.
$\Hinf$ norm evaluation involves repeated solution of several large-scale Riccati equations of order $n+r$.   Solving even a single Riccati equation, let alone several, will be a formidable task when $n$ is on the scale of  tens of thousands or larger, the range of system dimension of interest here.    We describe below an effective strategy to circumvent this difficulty. 

The optimization problem of Step 2 can be  rewritten (from Theorem \ref{thm:sameprojection}) as 
$$
\min \limits_{d_r \in \reals} \left \| H(s)-H_r^{0}(s)-\frac{d_r (G_1(s)-1)(G_2(s)-1)}{1-d_r G_3(s)} \Hinfty
$$
where $H_r^{0}(s)$ is a reduced-model obtained by \IRKA~in Step 1 of \textsf{{Algorithm IHA}}. 

 If one can find a reduced-order  approximation, $F_k(s)$,  to the error system, $F(s)=H(s)-H_r^{0}(s)$, having modest fidelity and order $k \ll n$, 
 then an associated optimal $d_r$-term could be efficiently calculated by solving the (comparatively) low order optimization problem 
\begin{equation} \label{eqn:approxHinf}
\min \limits_{d_r \in \reals} \left \|F_k(s)-\frac{d_r (G_1(s)-1)(G_2(s)-1)}{1-d_r G_3(s)} \Hinfty
\end{equation}
Provided  $k \ll n$, the cost of solving  (\ref{eqn:approxHinf}) will be negligible compared to that of original problem.  Of course, whatever advantage this strategy may bring  could be nullified if the cost of obtaining $F_k(s)$ is significant.  By using a Loewner matrix approach developed by Mayo and Antoulas  \cite{paramD} and described briefly below, we are able to obtain $F_k(s)$ at negligible cost relative to the computational demands already incurred in Step 1.  We reuse information obtained in the course of \IRKA~in Step 1 to obtain, for negligible additional effort, a reduced error model $F_k(s)$ having modest fidelity, adequate for the demands of Step 2. 


Suppose we have evaluated  the error system, $F(s)$, and derivative, $F'(s)$, on
a set of distinct points $\{s_1,\ s_2,\ \ldots,\ s_\ell\}\subset \IC$. We will construct  from this data a reduced order surrogate, 
$F_k(s) = \widehat{\boldc}_k^T(s \widehat{\boldE}_k-\widehat{\boldA}_k)^{-1}\widehat{\boldb}_k$ so that 
$$
F(s_i) = F_k(s_i) ~~{\rm and} ~~F'(s_i) = F'_k(s_i) ~~{\rm for} ~~i = 1,2,\ldots,\ell.
$$

The Loewner matrix approach as developed by Mayo and Antoulas  \cite{paramD} permits  ``data-driven" model reduction; 
one need not have access to state-space matrices determining a realization of the full order system.  
Only ``response measurements" are used, that is,  transfer function evaluations. 
Reduced-order models will be constructed directly that interpolate this ``measured data".  

Define  matrices $\IL \in \IC^{\ell \times \ell}$ and 
$\IM \in \IC^{\ell \times \ell}$ as
\begin{equation} \label{defineL}
\left(\IL \right)_{i,j} :=
\left\{  \begin{array}{ll}
\displaystyle 
  \frac{ F(s_i)-F(s_j)}{s_i-s_j}  & {\rm if}~~ i \neq j \\ \\
  F'(s_i) & {\rm if}~~ i = j 
\end{array} \right.
~~~
\left(\IM\right)_{i,j} := 
\left\{  \begin{array}{ll}
\displaystyle 
  \frac{s_i F(s_i)-s_jF(s_j)}{s_i-s_j}  &  {\rm if}~~i \neq j \\ \\
 \left. [sF(s)]'\right|_{s=s_i} & {\rm if}~~i= j 
\end{array} \right.
\end{equation}
$\IL$ is the \emph{Loewner matrix} associated with interpolation points $s_1,s_2,\ldots,s_\ell$ and the dynamical system $F(s)$; $\IM$ is the corresponding \emph{shifted Loewner matrix} (see \cite{paramD} for details). 
 Once $\IL$  and $\IM$ are constructed, assume that the interpolation data satisfy the following assumption:
\begin{equation} \label{LoewRankAssump}
{\rm rank}\,(s_i\IL-\IM )=
{\rm rank}[\IL~~\IM]={\rm rank} \left[ \begin{array}{c}\IL \\ \IM \end{array}\right]
\end{equation}
 for $i=1,2,\ldots, \ell$.  For SISO systems, this assumption holds whenever an interpolant of order $r={\rm rank}\,(s_i\IL-\IM )$ exists \cite{paramD}.  For MIMO systems, this assumption is also generically valid, but details are more involved; the interested reader should see Lemma 5.4 of \cite{paramD}.

In light of (\ref{LoewRankAssump}), a rational Hermite interpolant  is constructed by first choosing $k$ so that
 $ {\rm rank} \left[ \begin{array}{c}\IL \\ \IM \end{array}\right]\geq k$.
Then for some choice of $1\leq i \leq \ell$, compute
$
s_i\IL-\IM =\boldY\boldTheta \boldX^* ,
$
 the SVD of $s_i\IL-\IM$. 
Let $\boldY_k\in \IC^{\ell \times k}$ and $\boldX_k\in \IC^{\ell\times k }$ denote the leading $k$ columns of $\boldY$ and $\boldX$, respectively (associated with a truncated SVD of order $k$). 
Let $\boldZ = [F(s_1),~F(s_2)~,\ldots,~F(s_\ell)]^T$ and define 
 $$
 \widehat{\boldE}_k= -\boldY_k^*\IL \boldX_k,~~ 
 \widehat{\boldA}_k= -\boldY_k^* \IM \boldX_k, ~~
\widehat{\boldb}_k = \boldY_k^* \boldZ,~~\widehat{\boldc}_k = \boldZ^T \boldX_k,
$$
  and $F_k(s) =  \widehat{\boldc}_k^T(s \widehat{\boldE}_k-\widehat{\boldA}_k)^{-1}\widehat{\boldb}_k $.
$k$ may be considered as a truncation index here.   Depending on whether $k$ is chosen so that $k= {\rm rank}\,(s_i\IL-\IM )$ or $k< {\rm rank}\,(s_i\IL-\IM )$,
$F_k(s)$ will then be either an exact or an approximate interpolant,  respectively.
In practice, one should choose $k$ no larger than the numerical rank of 
$s_i\IL-\IM$, which can be determined by a Singular Value Decomposition (SVD) of $s_i\IL-\IM$.  
See \cite{paramD} for a full development of these ideas. 


Observe that until convergence occurs within Step 1 of \textsf{\small{Algorithm}}~\HinfIRKA, every cycle of Step 3 in \textsf{\small{Algorithm}}~\IRKA~will provide a sampling of $H(s)$ and $H'(s)$ at $r$ interpolation points -- generally a different set of interpolation points in each cycle. (Step 3 of  \textsf{\small{Algorithm}}~\IRKA~constructs a Hermite interpolant on a set of $r$ interpolation points that is cyclically adjusted.)  Suppose that \textsf{\small{Algorithm}}~\IRKA~takes $q$ steps to convergence.  When Step 1 of \textsf{\small{Algorithm}}~\HinfIRKA~concludes, we will have had $H(s)$ and $H'(s)$ sampled at a total of $\ell = q \times r$ interpolation points. We collect these interpolation points and transfer function evaluations throughout \IRKA. Once \IRKA~converges, (yielding an $\mathcal{H}_2$-optimal model, $H_r^0(s)$), we evaluate $H_r^0$ and $H_r^{0\,\prime}$  at these $\ell$ points as well. Since the order of $H_r^0(s)$ is $r$, the cost of these function evaluations is negligible.  After the completion of Step 1 of \textsf{\small{Algorithm}}~\HinfIRKA, we have 
an $\ell$-fold sampling of both $F(s) = H(s) - H_r^0(s)$ and $F'(s)$ with virtually no additional computational cost beyond what was needed for Step 1 itself. 
Then, we simply apply the Loewner matrix approach described above to construct $F_k(s)$.
The choice of $k$ will be clarified via numerical examples in \S\ref{sec:examples}. 

\paragraph{Numerical Cost of \HinfIRKA.}
Note that once Step 1 of \textsf{\small{Algorithm}}~\HinfIRKA~is completed, Step 2 is not computationally intense.  The cost is dominated typically by an $\ell \times \ell$ SVD computation
with $\ell = q \times r$.  Since \IRKA~typically converges rather quickly (especially so in the SISO case focused on here), $\ell$ is generally modest in size.
In all of our numerical examples, we have never needed to compute an SVD of size larger than $200 \times 200$; a trivial computation. Moreover, in all of our numerical examples $k$ never exceeded $33$; making the solution of the optimization problem in Step 2 quite cheap. The overall cost of \HinfIRKA\ is only marginally more than that of \IRKA\ and is dominated by the same sequence of sparse linear solves.


\paragraph{Stability of the reduced model:}
The asymptotic stability of $H_r$ in Step 2 may be enforced by adding  a penalty function to the cost function penalizing values of $d_r$ that yield systems having poles too close to the imaginary axis.  The optimization algorithm would then automatically reject $d_r$  terms that cause unstable eigenvalues in the pencil $s\boldEr-\boldAr - d_r\bolde \bolde^T$.  
Alternatively, one may simply calculate and check the eigenvalues of this $r\times r$ pencil to determine whether to accept a $d_r$ on the basis of stability. In our numerical examples, we simply reject $d_r$ values that cause unstable reduced models by setting the corresponding the function value to $\infty$. We always obtained a stable reduced-model as a result but there are better, more effective numerical strategies  to perform this task.  For example, a logarithmic barrier function that takes the real part of the pole closest to the imaginary axis as its argument. These numerical issues will be studied in a separate work where we extend the method to the MIMO as discussed next.

\paragraph{Application to MIMO systems:}
\textsf{\small{Algorithm}}~\HinfIRKA, together with the results of Section 2.2, can be easily generalized to the MIMO case.  In the MIMO case, \IRKA\ enforces bitangential interpolation conditions at the reflection of the reduced-order poles across the imaginary axis.  A complete account of $\mathcal{H}_2$ optimal model reduction and \IRKA\ for MIMO systems can be found in \cite{Ant2010imr}.  These interpolation conditions can be enforced while varying the D-term of the reduced-order system.  See Theorem 3 in \cite{Struct} for a complete description of how to construct the bitangential interpolant while varying the D-term; thus the theory in this paper directly generalizes to the MIMO case.  The optimization step involving the matrix $\boldsymbol{D}\in \reals^{p\times m}$, however, is more involved than its scalar counterpart for SISO systems and a robust numerical implementation  is the subject of ongoing research. 

\section{Numerical Experiments} \label{sec:examples}

We illustrate here the performance of  \HinfIRKA~on various benchmark models for model reduction and compare its performance with that of  balanced truncation (\BT)~ optimal Hankel norm approximation (\HNA), and Iterative Rational Krylov Algorithm (\IRKA). We note that generic \BT~yields $d_r=0$, a null feed-forward term. Therefore, in order to present  a fair
representation for \BT, we add a final step that varies the $d_r$ term in the \BT~model as well; we use
$d_r$  that gives the minimum $\Hinf$ norm.   We refer to \BT~with this optimally chosen $d_r$ term as ``modified balanced truncation'' (or \MBT).

\subsection{PEEC Model}
\begin{figure}[hhhht]
\centerline{\hbox{
\epsfysize=8cm
\epsfxsize=10cm
\epsffile{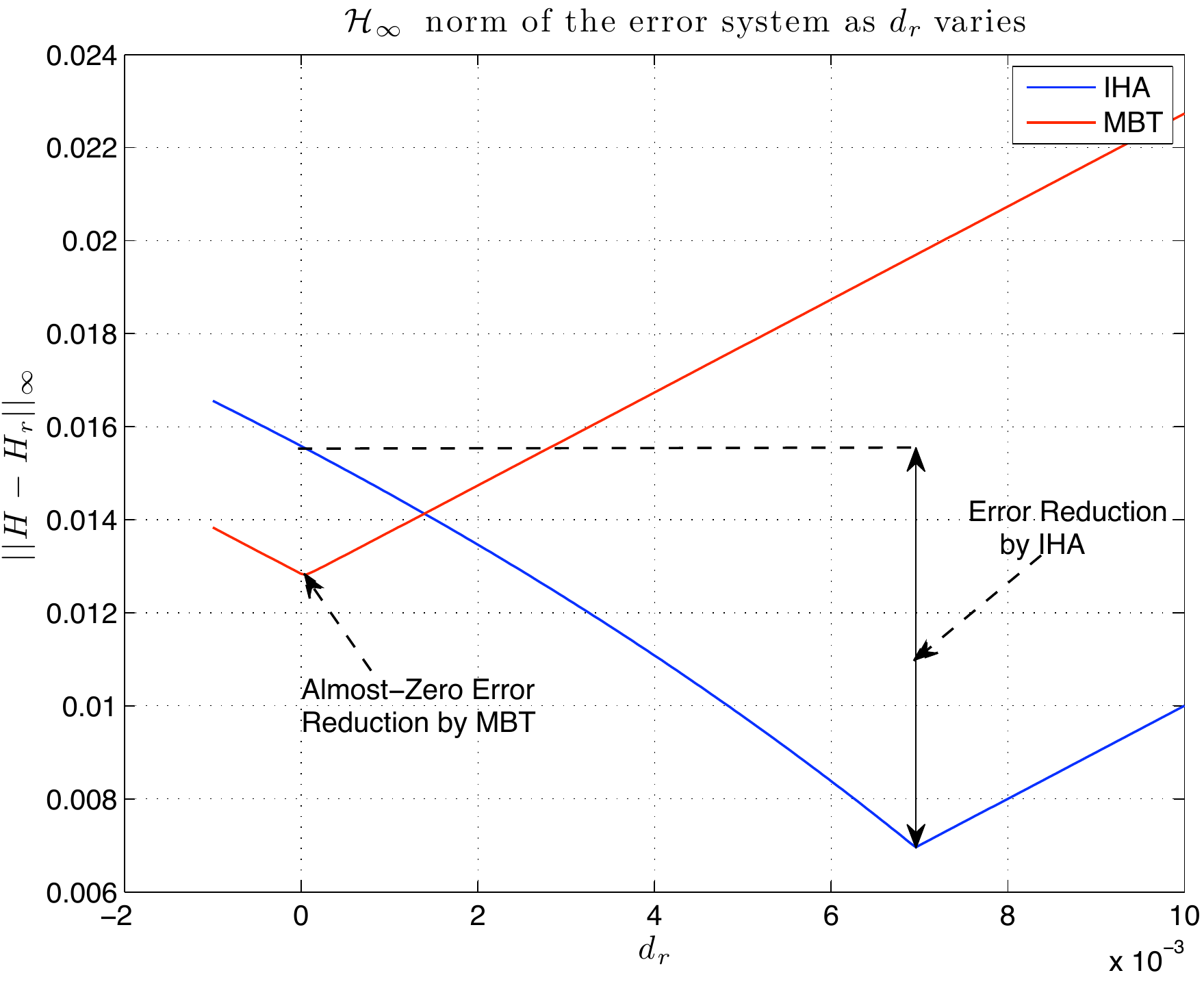}
}
}
\caption{Comparison of the $\Hinf$ Error as the \dr-term varies for \HinfIRKA~and \MBT}\label{dr_term_error_evolution}
\end{figure}

The full-order system is the spiral inductor system PEEC model \cite{chahlaoui2005benchmark} 
of order $n=1434$. The system is state-space-symmetric (\SSS), i.e.
the transfer function 
$ \Hs= \boldc^T(s\boldE-\boldA)^{-1}\boldb$ satisfies  $\boldE=\boldE^T >0$, $\boldA=\boldA^T$ and $\boldc=\boldb$. \SSS~systems appear in many important applications such as in the analysis of RC circuits, and has been the subject of several model reduction papers; e.g, see \cite{zippreserve,Zipchar,reis2011lbp,bai2006hmr}.
We first illustrate the  effect of the $d_r$-term modification for 
\MBT~and for Step 2 of \textsf{\small{Algorithm}}~\HinfIRKA. In \MBT, once the initial \BT~phase is completed, we vary $d_r$ and measure the resulting changes to the $\Hinf$ error. 
For \HinfIRKA, once the \IRKA~phase in Step 1 is completed, we vary the $d_r$ term (inducing corresponding changes
 to $\mathbf{A}_r$, $\mathbf{b}_r$ and $\mathbf{c}_r$ as in Theorem \ref{thm:sameprojection}) and once again
measure the resulting changes to the $\Hinf$ error. Note that calculation of $\Hinf$ error is 
for illustration  only and is not used by \HinfIRKA~to compute the optimal $d_r$-term. 
Results are shown in 
Figure \ref{dr_term_error_evolution} for $r=2$.  As the figure shows,  there is essentially no improvement in \MBT~that comes from adjusting $\dr$, 
however for \HinfIRKA\ the error is reduced by a factor larger than two.
Even though at  the starting point, $d_r=0$, the $\Htwo$-optimal approximation has higher $\Hinf$ error  than does \BT\ at $d_r=0$,  \HinfIRKA\ 
is able to reduce the $\Hinf$ error to a value significantly lower than that for \MBT\ through $d_r$-term optimization, at virtually negligible computational cost. 
In Figure \ref{dr_term_error_evolution}, the point $d_r=0$ on the curve for \HinfIRKA\ gives the value of the $\Hinf$ error produced by \IRKA. Note that the $\Hinf$ error for \HinfIRKA\ is less than half of that for \IRKA. This behavior is common to all the numerical examples that follow.

\begin{figure}[hhhht]
\centerline{\hbox{
\epsfysize=8cm
\epsfxsize=10cm
\epsffile{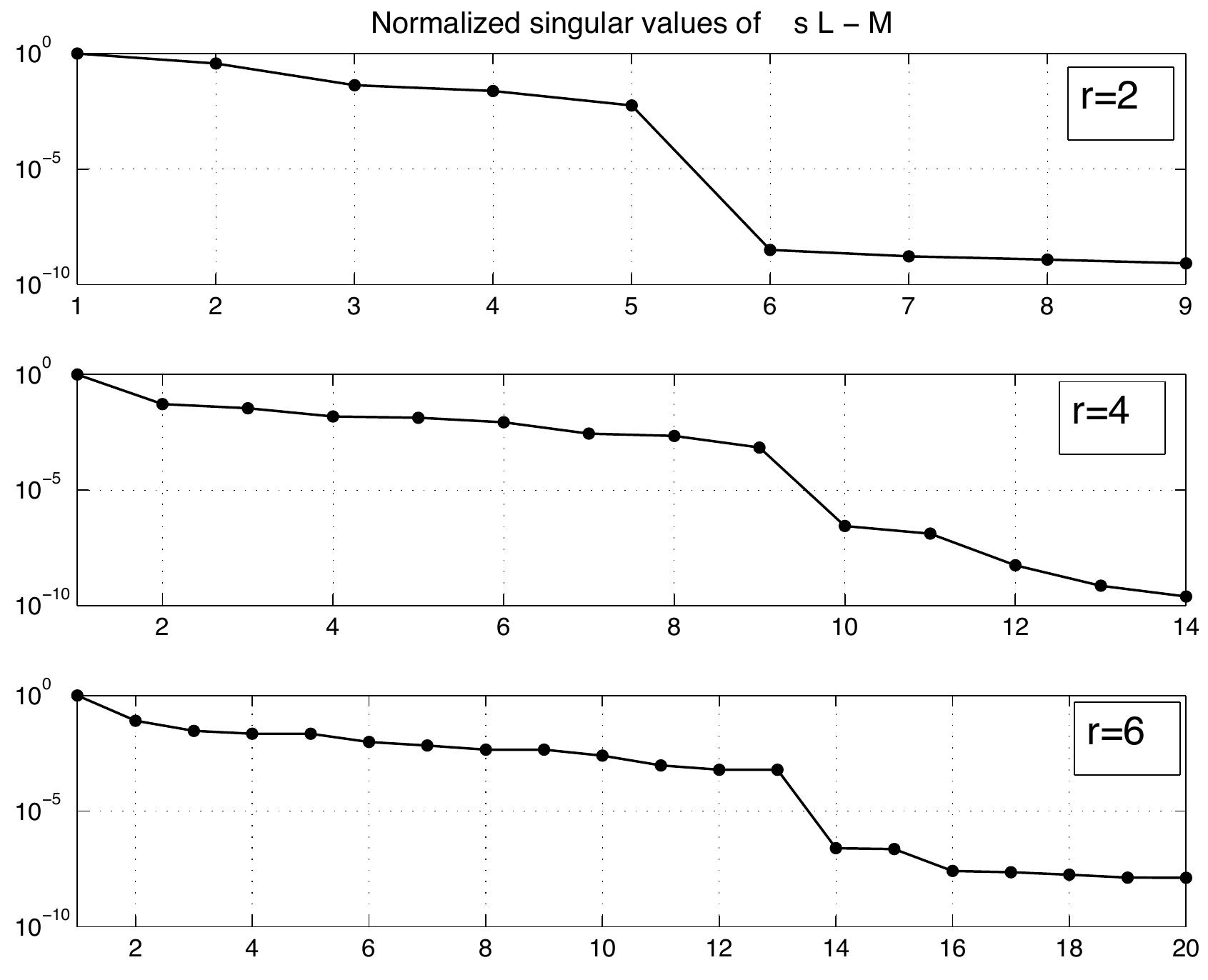}
}
}
\caption{The decay of the singular values of $s_i\IL-\IM$ for PEEC Model}\label{fig:peec_loewner_decay}
\end{figure}

Before presenting the comparisons between \HinfIRKA, \MBT, and \HNA, we illustrate the efficiency of the methodology outlined in \S \ref{loewner} in solving the optimization problem in \textsf{\small{Algorithm}}~\HinfIRKA; in other words in finding the optimal $d_r$-term in Figure \ref{dr_term_error_evolution}. For the three $r$ values $r=2,4,6$, we implement \HinfIRKA~both by exactly solving 
Step 2 and by the method of \S \ref{loewner}. Table \ref{table:peec:loewner-vs-exact} tabulates the results where the resulting optimal $d_r$ values and the $\Hinf$ error norms for both methods together with the order-$k$ used in  method of \S \ref{loewner} are listed. As the table clearly illustrates that the Loewner matrix approach yields $d_r$ terms and the $\Hinf$ error norms which are very close to true-optimal values of the underlying optimization problem in Step 2 of \textsf{\small{Algorithm}}~\HinfIRKA. More importantly, this is achieved with negligible computational cost where the function evaluations are the $\Hinf$ norm computations for an order $k$ system only as opposed to order $n+r$. 
\begin{table}[htb]
\caption{Solution of the optimization problem in Step 2}
\center
\begin{tabular}{c|cc|ccc}
 \multicolumn{3}{c|}{Exact} &\multicolumn{2}{c}{Loewner}  \\ \hline
$r$&$d_r^\star$&$\| H - H_r^\star \|_{\Hinf}$&$d_r^\star$&$\| H - H_r^\star \|_{\Hinf}$ & $k$ \\
\hline
$2$& $6.9577 \times 10^{-3}$  & $4.4522 \times 10^{-3} $ & $6.9659 \times 10^{-3}$ & $4.4574 \times 10^{-3}$ & $5$\\  
$4$& $1.0041\times 10^{-4}$  & $8.6577 \times 10^{-5}$ & $1.0076\times 10^{-4}$ & $8.7114\times 10^{-5}$ & $9$\\
$6$ & $2.7795 \times 10^{-6}$ & $4.4771 \times 10^{-6}$&  $2.7804 \times 10^{-6}$& $4.4857 \times 10^{-6}$ & $13$\\
\end{tabular}
\label{table:peec:loewner-vs-exact}
\end{table}
One pattern we have observed via several \SSS~models is that $k=2r+1$ is a natural choice. This pattern repeated itself for every \SSS~example we have tried. There was a clear cut-off point in the singular values of the matrix $s_i \IL - \IM$ at $(2r+1)^{\rm th}$ singular value. The decay of these singular for all three $r$ values are shown in Figure \ref{fig:peec_loewner_decay} supporting the $k=2r+1$ choice. 

\begin{figure}[hhhht]
\centerline{\hbox{
\epsfysize=6cm
\epsfxsize=8cm
\epsffile{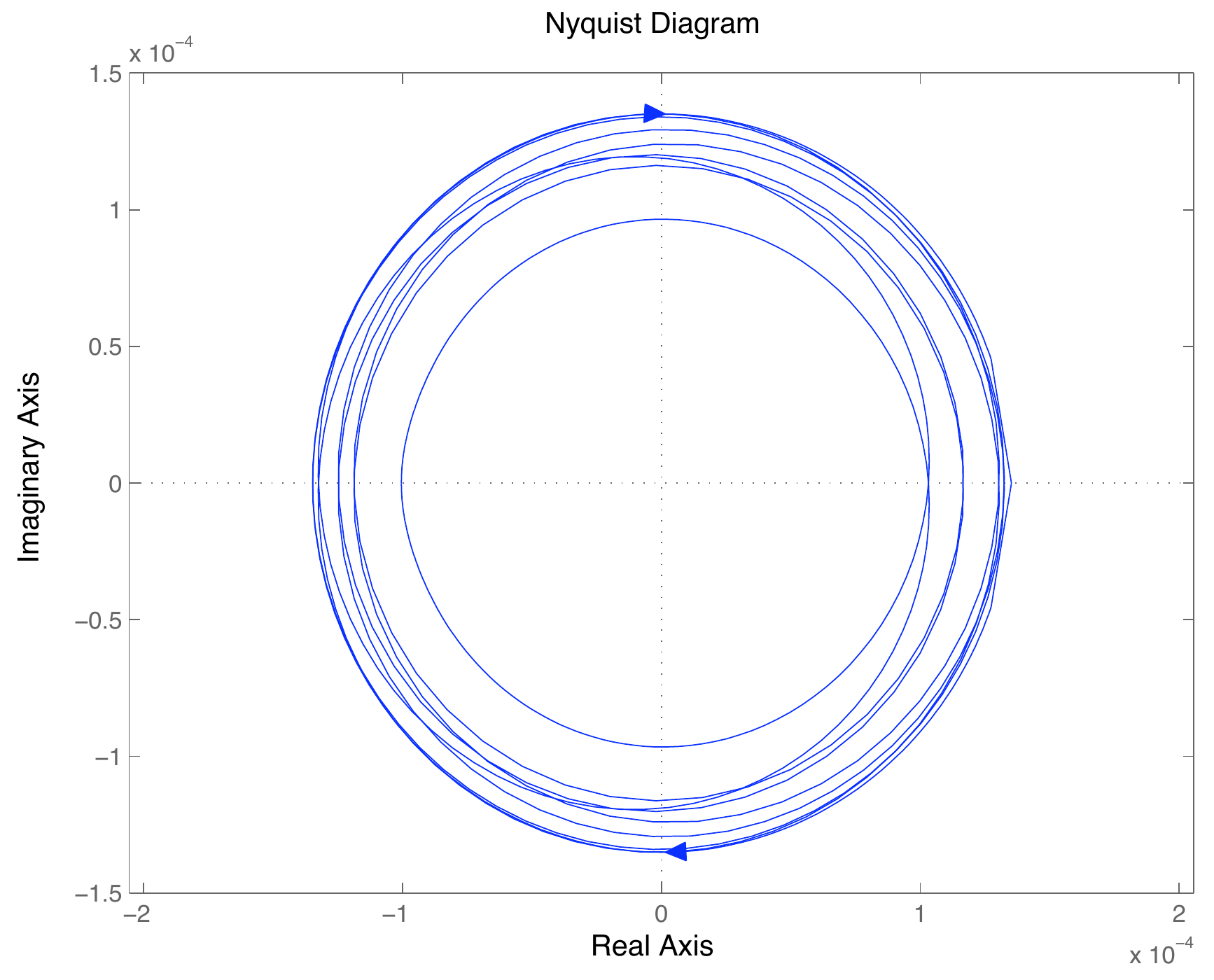}
}
}
\caption{The Nyquist plot of the error system $H-H_r^\star$ for PEEC Model}\label{Nyquistpeec}
\end{figure}

Next, to compare \HinfIRKA, \MBT, and \HNA,
we reduce the order of the system to $r=2,4,6$. 
The resulting relative $\Hinf$ error values together with 
the lower  bound (i.e. $\sigma_{r+1}/\| H \|_{\Hinf}$ where $\sigma_{r+1}$ denotes the $(r+1)^{\rm th}$ Hankel singular value of $H(s)$)  are listed in Table \ref{peec_table}. The lowest error value for each $r$
is shown in bold font. 
For every $r$, \HinfIRKA~outperforms \MBT~by almost a factor of two. Also it performs very close to \HNA, indeed outperforms it for $r=4$. This shows the strength of the proposed method. Without solving any Lyapunov equations and without the need for any large-scale $\Hinf$ norm computation, our methods consistently outperforms \BT~by a significant amount and performs as nearly as and sometimes better than 
\HNA. Note that $\Hinf$ error values for \ \HinfIRKA~is very close to the lower bound given by $\sigma_{r+1}/\| H \|_{\Hinf}$. We can relate this to Theorem \ref{thm:trefsufficient}: 
  for each order of approximation shown in Table \ref{peec_table}, the reduced-model due to  \HinfIRKA~results in exactly
$2r+1$ interpolation points in the right-half plane and a 
nearly circular error curve as illustrated in
 Figure \ref{Nyquistpeec} .  
 $2r$ of these zeros result from \IRKA; indeed these are $r$ distinct zeroes with multiplicity $2$. Then, the 
 $(2r+1)^{\rm th}$ zero is obtained by
introducing the $d_r$ term. 
 For example, for the case $r=6$, computing the optimal \dr-term in turn placed an additional zero at the point $6.60 \times 10^{8}$.  
\begin{table}[htb]
\caption{Relative $\Hinf$ error norms  for the PEEC Model}
\center
\begin{tabular}{cccc|c}
$r$&\HinfIRKA&\MBT&\HNA& Lower bound \\
\hline
$2$&$4.45\times 10^{-3}$ &$8.21\times 10^{-3}$&$\mathbf{3.95 \times 10^{-3}}$&$3.72 \times 10^{-3}$ \\
$4$&$\mathbf{8.66 \times 10^{-5}}$&$2.78 \times 10^{-4}$&$1.24\times 10^{-4}$&$7.79\times 10^{-5}$\\
$6$ &$4.48 \times 10^{-6}$&$8.16\times 10^{-6}$&$\mathbf{3.41\times10^{-6}}$&$3.15\times 10^{-6}$
\end{tabular}
\label{peec_table}
\end{table}

Note that the findings in this \SSS~example are common to all other \SSS~examples we have tried. The $\Hinf$ error due to \HinfIRKA~has  always been around half of \MBT, and $k=2r+1$ has been a clear cut-off point for the Loewner matrix approach. Due to space limitations, we omit these examples. For 
some other \SSS~examples, we refer the reader to \cite{flagg2009iba}.

\subsection{CD Player Model}
\begin{figure}[hhhht]
\centerline{\hbox{
\epsfysize=8cm
\epsfxsize=10cm
\epsffile{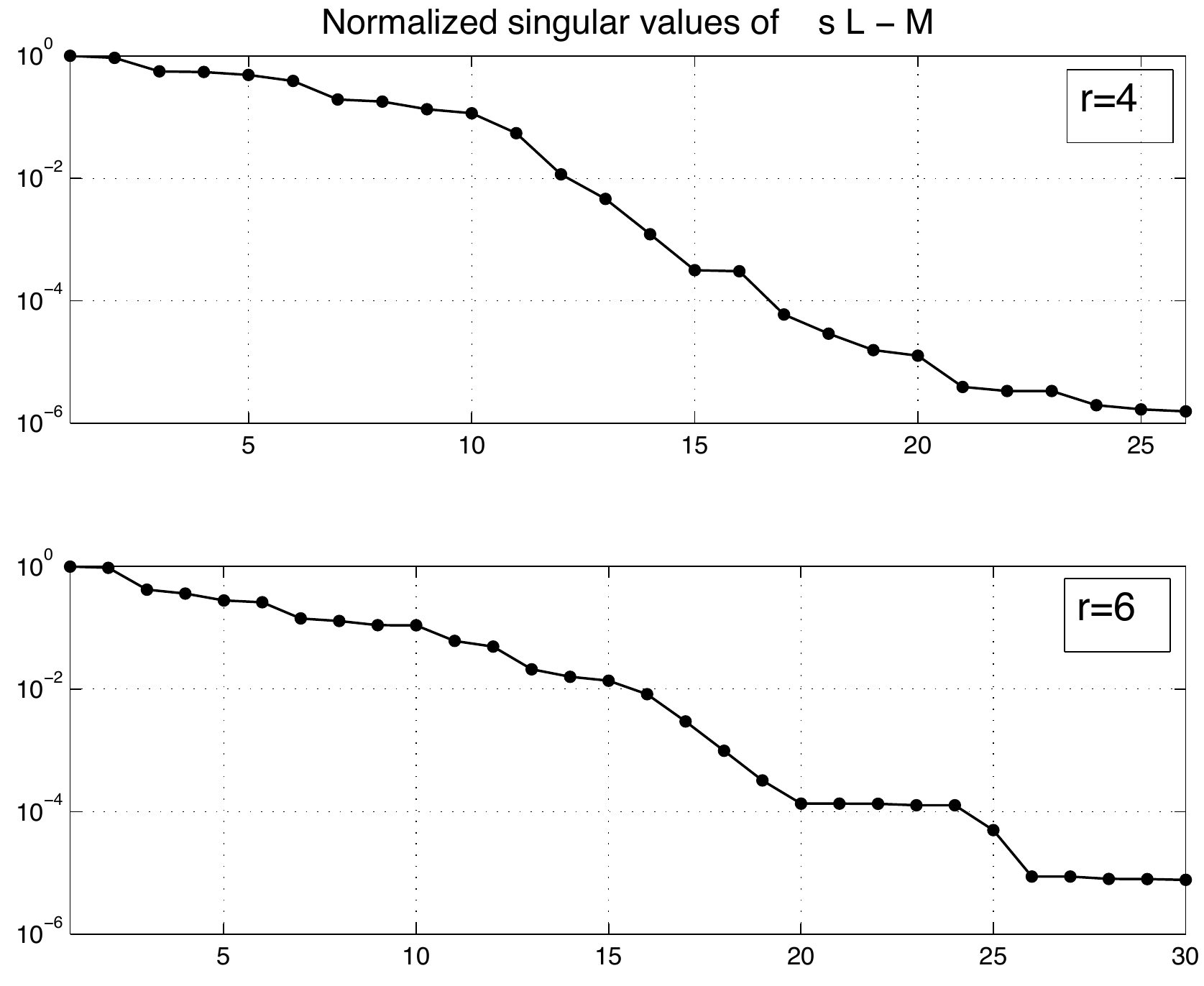}
}
}
\caption{The decay of the singular values of $s_iL-M$ for the CD Player Model}\label{fig:cd_loewner_decay}
\end{figure}

Next, we demonstrate the proposed method on the CD player model of order $n=120$. For details on this model, see \cite{antoulas2001asurvey,chahlaoui2005benchmark}.  We reduce the order to $r=2,4,6,8,10$ using \HinfIRKA; stopping at $r=10$ as the relative error fell below $10^{-3}$.  As for the previous example,  we compare the effect of 
solving the optimization problem in Step 2 of  \textsf{\small{Algorithm}}~\HinfIRKA\ exactly versus by the method of  \S \ref{loewner}.
The results are listed in Table \ref{table:cd:loewner-vs-exact}. As for the previous case, 
 the Loewner matrix approach yields $\Hinf$ error norms which are very close to the true optimal values of the underlying optimization problem. Once again, the function evaluations are much simpler as $k$ never exceeds $33$. In this example, we chose the value of $k$ as the normalized singular values of $s_i\IL-\IM$ drops below 
 the tolerance of $10^{-5}$. Figure \ref{fig:cd_loewner_decay} shows this decay behavior for $r=4$ and $r=6$. 
 \begin{table}[htb]
\caption{Solution of the optimization problem in Step 2}
\center
\begin{tabular}{c|lc|lcc}
 \multicolumn{3}{c|}{Exact} &\multicolumn{2}{c}{Loewner}  \\ \hline
$r$&$~~~~~~d_r^\star$&$\| H - H_r^\star \|_{\Hinf}$&$~~~~~~d_r^\star$&$\| H - H_r^\star \|_{\Hinf}$ & $k$ \\
\hline
$2$   & $-3.6728$                        & $3.6597 \times 10^{-1}$  & $-3.6545$                        & $3.6604 \times 10^{-1}$   & $2$\\  
$4$   & $\phantom{-}3.5019\times 10^{-1}$   & $2.1318 \times 10^{-2}$  & $\phantom{-}2.4819\times 10^{-1}$    & $2.1422 \times 10^{-2}$    & $20$\\
$6$   & $\phantom{-}2.9538 \times 10^{-1}$  & $1.0155 \times 10^{-2}$  &  $\phantom{-}2.0763 \times 10^{-1}$  & $1.0426 \times 10^{-2}$   & $25$\\
$8$   & $\phantom{-}1.3888 \times 10^{-1}$  & $4.8357 \times 10^{-3}$  &  $\phantom{-}1.3625 \times 10^{-1}$  & $4.8526 \times 10^{-3}$   & $32$\\
$10$ & $-3.5750 \times 10^{-2}$ & $8.5384\times 10^{-4}$  &  $-3.2438 \times 10^{-2}$ & $8.9952\times 10^{-4}$   & $33$\\
\end{tabular}
\label{table:cd:loewner-vs-exact}
\end{table}

Next, we compare  \HinfIRKA~with \MBT~and \HNA.
The results are illustrated in Table \ref{CDtable}  where the minumum 
value for each $r$ is shown in bold font. 
Note that for every $r$ value, the proposed approach outperforms \MBT. Moreover, except for the $r=2$ and $r=4$ cases, \HinfIRKA~outperforms 
\HNA~as well. 

\begin{figure}[hht]
\centerline{\hbox{
\epsfysize=8cm
\epsfxsize=10cm
\epsffile{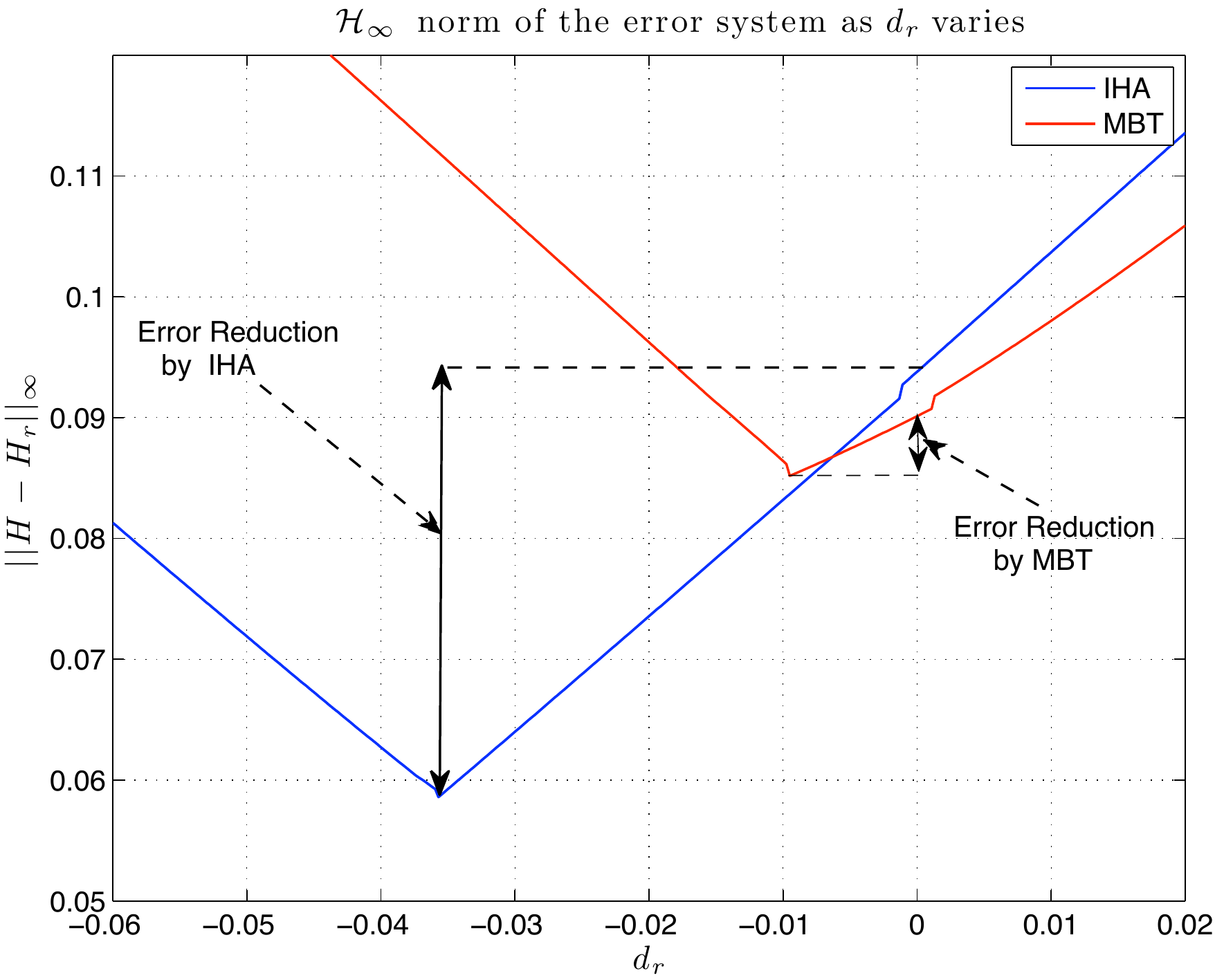}
}
}
\caption{$\Hinf$ error as a function of the $\dr$-term in \MBT~and \HinfIRKA}\label{cd_dr_gain}
\end{figure}

Before moving to the next example, we illustrate the effect of the $d_r$-term modification in our proposed method as opposed to \MBT. 
In Figure \ref{cd_dr_gain}, we show how the {\it absolute} $\Hinf$-error changes both in  \HinfIRKA~and in \MBT~as we vary the $d_r$-term for $r=10$. In both cases, $d_r=0$ is the starting point. While the $\Hinf$-error
reduces marginally from $0.0905$ to $0.0852$ in \MBT--\,only a $5.86\%$ reduction, the gain is much more significant in  \HinfIRKA~where we reduce the $\Hinf$-error from $0.0938$ down to $0.0585$, a significant $37.67\%$ reduction in the $\Hinf$ error. Even though the $\Htwo$-optimal approximation has a larger $\Hinf$ error at the starting point, $d_r=0$, than does the \BT\ approximation, 
 the $d_r$-term optimization in \HinfIRKA\ will produce a final $\Hinf$ error that is significantly lower than that produced by the corresponding $d_r$-term optimization in \MBT. 
\begin{table}[hhhht]
\caption{ Relative \Hinfspace-norm Error norms for the CD Player Model}\label{CDtable}
\center
\begin{tabular}{c|ccc}
$r$ & \HinfIRKA & \MBT &\HNA \\
\hline
$2$&$3.66 \times 10^{-1}$&$3.68\times 10^{-1}$&$\mathbf{3.35\times 10^{-1}}$\\
$4$&$2.14\times 10^{-2}$&$2.25\times 10^{-2}$&$\mathbf{2.00\times 10^{-2}}$\\
$6$&$\mathbf{1.04\times10^{-2}}$&$1.19\times 10^{-2}$&$1.23\times 10^{-2}$\\
$8$&$\mathbf{4.85\times10^{-3}}$&$6.40\times10^{-3}$&$5.99\times10^{-3}$\\
$10$&$\mathbf{8.99\times10^{-4}}$&$1.24\times10^{-3}$&$1.08\times10^{-3}$
\end{tabular}
\end{table}

\subsection{Heat Model} 
\begin{figure}[hhhht]
\centerline{\hbox{
\epsfysize=6cm
\epsfxsize=10cm
\epsffile{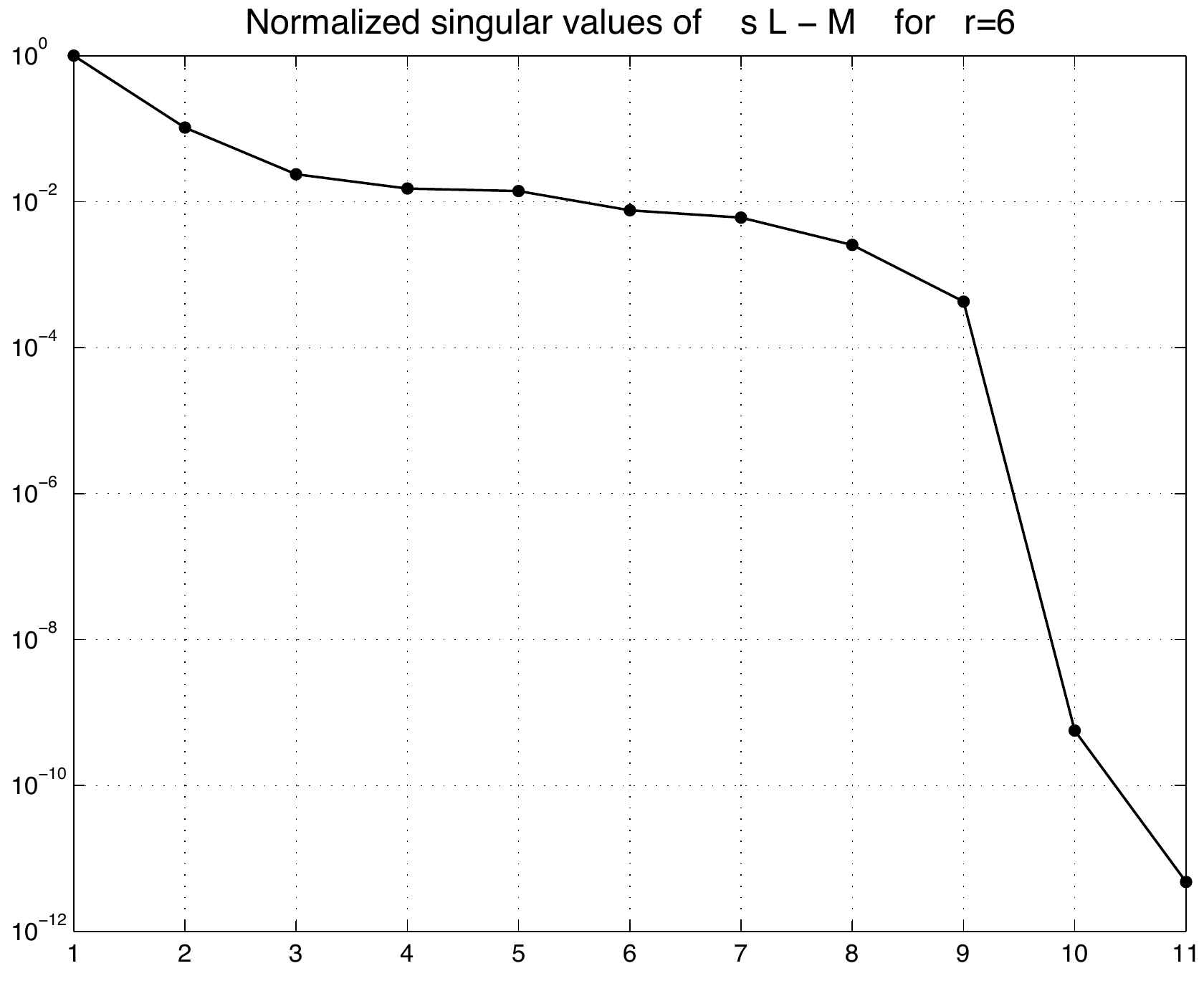}
}
}
\caption{The decay of the singular values of $s_i\IL-\IM$ for the Heat Model}\label{fig:heat_loewner_decay}
\end{figure}

The full-order model  is a plate with two heat sources and two 
points of measurements, and described by the heat equation as explained in \cite{antoulas2001asurvey,gugercinprojection}. A model of 
order $n=197$ is obtained by spatial discretization. We choose a SISO subsystem corresponding to 
the first input and first output.  Using \HinfIRKA, we reduce the order to $r=2,4$ and $r=6$. 
As in the previous examples, we first tabulate, in Table \ref{table:heat:loewner-vs-exact}, the results for solving the
Step 2 of  \textsf{\small{Algorithm}}~\HinfIRKA\ exactly and by the method of  \S \ref{loewner}. The conclusion is the same as before: 
 The Loewner matrix approach yields $\Hinf$ error norms and optimal $d_r$ values which are very close to the true optimal values of the underlying optimization problem in Step 2 of \textsf{\small{Algorithm}}~\HinfIRKA, indeed exact to the fifth digit for the $r=6$ case. 
 \begin{table}[htb]
\caption{Solution of the optimization problem in Step 2}
\center
\begin{tabular}{c|lc|lcc}
 \multicolumn{3}{c|}{Exact} &\multicolumn{2}{c}{Loewner}  \\ \hline
$r$&$~~~~~~d_r^\star$&$\| H - H_r^\star \|_{\Hinf}$&$~~~~~~d_r^\star$&$\| H - H_r^\star \|_{\Hinf}$ & $k$ \\
\hline
$2$   & $-2.0694\times 10^{-1}$  & $1.0710 \times 10^{-2}$  & $-2.0700\times 10^{-1}$  & $1.0711 \times 10^{-2}$   & $7$\\  
$4$   & $\phantom{-}2.0875\times 10^{-2}$   & $8.9082 \times 10^{-4}$  & $\phantom{-}2.0813\times 10^{-2}$    & $8.9166 \times 10^{-4}$    & $9$\\
$6$   & $-5.6250 \times 10^{-4}$  & $2.3578 \times 10^{-5}$  &  $-5.6250 \times 10^{-1}$  & $2.3578 \times 10^{-5}$   & $13$
\end{tabular}
\label{table:heat:loewner-vs-exact}
\end{table}
The decay of the singular values of $s_i\IL-\IM$ is shown in Figure \ref{fig:heat_loewner_decay}, illustrating the choice of $k$. Only the $r=6$ case is presented; the other cases show the same pattern.

Results for comparison with \MBT~and \HNA~are shown in Table \ref{matchpointscompare}. Once more, the proposed method consistently yields better $\Hinf$ performance than \MBT.  Even though 
for $r=2$ the proposed method leads to smaller $\Hinf$ error, 
for $r=4,6$, \HNA~yields slightly better results.  Hence,
for this example as well, using interpolatory projections, we are able to beat \MBT~consistently and yield results comparable to or better than  that of \HNA. Indeed, this is satisfactory since, as stated in the introduction, for the large-scale settings we are interested in, implementing \HNA~will be a formidable task, if not impossible. 
\begin{table}[htb]
\caption{Relative $\Hinf$ error norms  for the Heat Model }
\label{matchpointscompare}
\center
\begin{tabular}{cccccc}
$r$&\HinfIRKA&\MBT&\HNA\\
\hline
$2$&$\mathbf{1.08\times 10^{-2}}$&$1.66\times 10^{-2}$&$1.11\times 10^{-2}$  \\
$4$&$8.92 \times 10^{-4}$& $1.68 \times 10^{-3}$&$\mathbf{8.47 \times 10^{-4}}$\\
$6$&$2.30\times 10^{-5}$&$4.61 \times 10^{-5}$ &$\mathbf{2.07 \times 10^{-5}}$
\end{tabular}
\end{table}

\begin{figure}[hht]
\centerline{\hbox{
\epsfysize=8cm
\epsfxsize=10cm
\epsffile{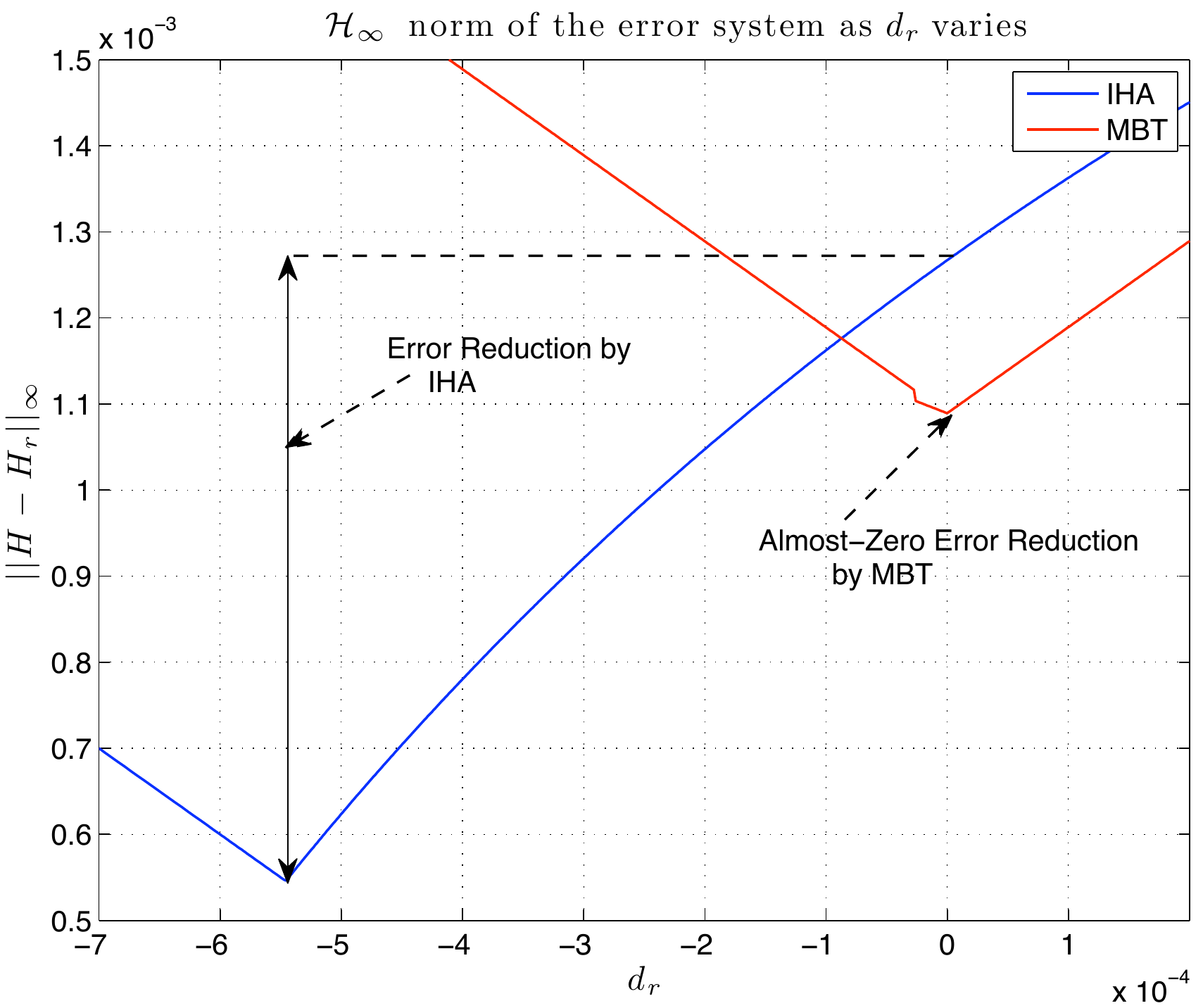}
}
}
\caption{$\Hinf$ error as a function of the $\dr$-term in \MBT~and \HinfIRKA}\label{heat_dr_gain}
\end{figure}

As done in the previous examples, 
we illustrate, in Figure \ref{heat_dr_gain}, the behavior of the {\it absolute} $\Hinf$-error while optimizing over the $d_r$-term in both \HinfIRKA~and \MBT. In this case, \MBT~almost gains nothing from the $d_r$-term modification, as the $\Hinf$-error is reduced from $1.0897\times 10^{-3}$ only to $1.0894\times 10^{-3}$, a marginal gain of $0.027\%$. 
On the other hand, \HinfIRKA~reduces the $\Hinf$-error from $1.26 \times 10^{-3}$ down to $5.46 \times 10^{-4}$, a reduction factor of $56.89\%$. Once again, this reduction in the error is achieved even when the initial $\Hinf$-error is bigger than that of \BT, providing a clear illustration of the effectiveness of the 
$d_r$-term optimization in the proposed method.

To illustrate the importance of using the \IRKA~points in Step 2 of \textsf{\small{Algorithm}}~\HinfIRKA, we use arbitrarily chosen interpolation points (within the bounds of the mirror spectrum of the full order \boldA) 
rather than the \IRKA~points. Then, we apply the same $d_r$-term modification as before. For $r=6$, for example, the resulting relative $\Hinf$-error is $1.72\times10^{-3}$,  two order of magnitudes higher than what we obtain using the \IRKA~points. This simple example  illustrates the advantage of initializing  Step 2 of \textsf{\small{Algorithm}}~\HinfIRKA~with points computed by \IRKA.  We want to emphasize that these arbitrary interpolation points are indeed used to initialize 
\IRKA. Hence, \IRKA~corrects these arbitrarily chosen points, producing interpolation points that are used to obtain an $\Hinf$-error norm two orders of magnitude smaller.

\subsection{A Heat Transfer Problem in the Cooling
 of Steel Profiles}

\begin{figure}[hhhht]
\centerline{\hbox{
\epsfysize=6cm
\epsfxsize=10cm
\epsffile{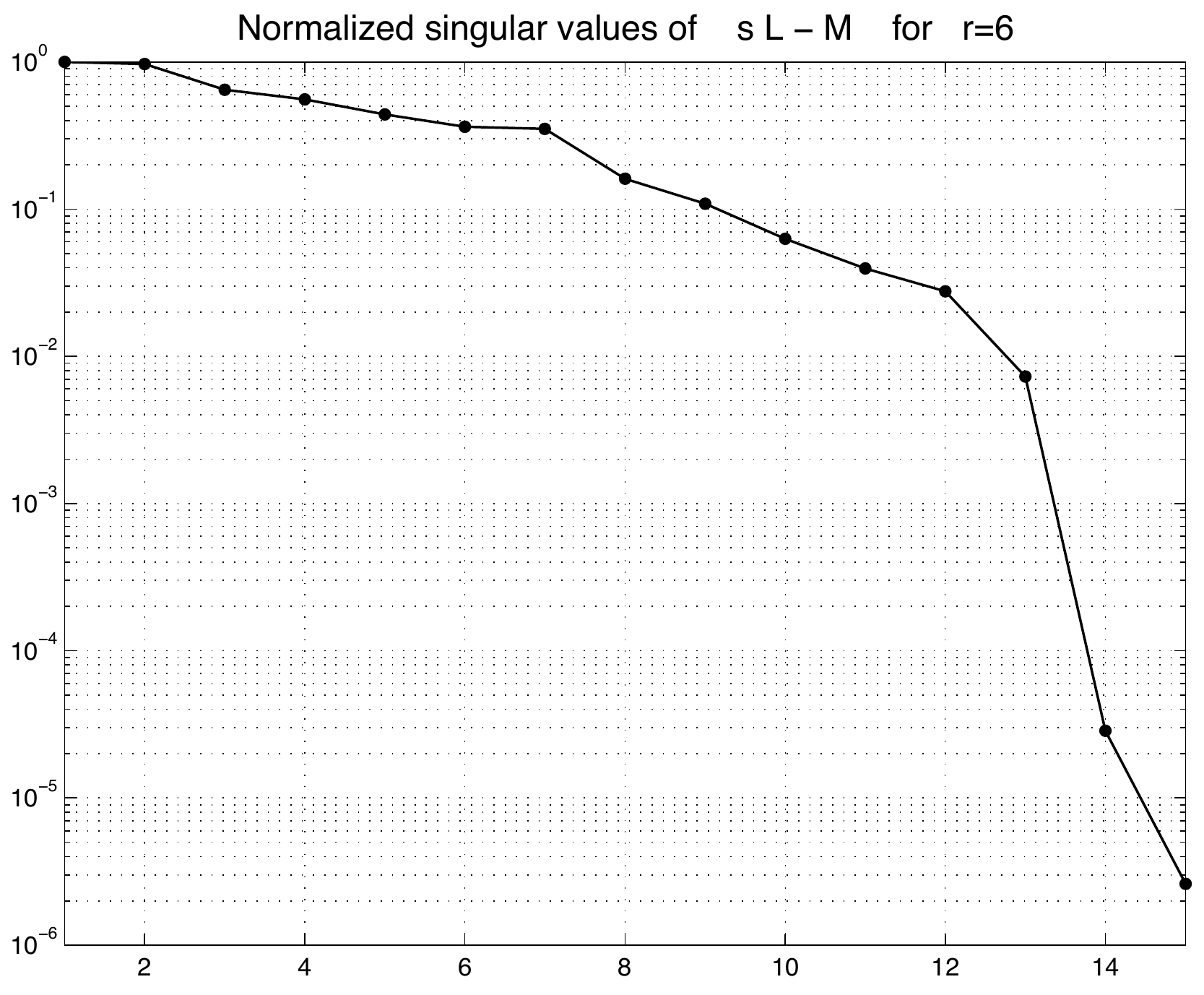}
}
}
\caption{The decay of the singular values of $s_i\IL-\IM$ for the Cooling Steel Model}\label{fig:irka_loewner_decay}
\end{figure}

This model describes a cooling process in
 a rolling mill which has been modeled as boundary control of a two dimensional heat equation. 
 The full order model has $n=79841$ with $7$ inputs and $6$ outputs. We consider a SISO subsystem corresponding to the sixth input and the second output.
For details regarding this model, see \cite{benner2004solving,benner2004ens}.
For such a large order, implementing \HNA~is not possible. \BT~can be implemented iteratively using ADI-type methods; however this requires state-of-the-art iterative Lyapunov solvers (two generalized Lyapunov equations of order $n=7984$1 need to be solved) and is not the focus of this paper. Hence, we concentrate only on the performance of the proposed method and compare it with \IRKA~to show the improvement by the $d_r$-term
optimzation. \IRKA~in Step 1 of \textsf{\small{Algorithm}}~\HinfIRKA~is implemented in Matlab using direct sparse linear solves. Once again, we consider the solution of the optimization problem in Step 2 of \textsf{\small{Algorithm}}~\HinfIRKA~using both approaches, i.e., directly solving the large-scale optimization problem versus using the method of \S\ref{loewner}. The results are shown in Table \ref{table:rail:loewner-vs-exact} and reveal the same pattern as before. Note that the direct solution of this scalar optimization problem requires several $\Hinf$ norm computation for a system of order $79841+r$. This is computationally intractable, so 
 the $\Hinf$ norms for the direct method are computed approximately by sampling along the imaginary axis at $500$ points, logarithmically spaced points between $10^{-8}$ and $10$.  
 This is not an issue for the Loewner matrix approach since  $\Hinf$-norm computations are done on surrogate error systems of order $k$; which does not exceed $13$ for this example.
 \begin{table}[htb]
\caption{Solution of the optimization problem in Step 2}
\center
\begin{tabular}{c|lc|lcc}
 \multicolumn{3}{c|}{Exact} &\multicolumn{2}{c}{Loewner}  \\ \hline
$r$&$~~~~~~d_r^\star$&$\| H - H_r^\star \|_{\Hinf}$&$~~~~~~d_r^\star$&$\| H - H_r^\star \|_{\Hinf}$ & $k$ \\
\hline
$2$   & $\phantom{-}1.0189\times 10^{-2}$  & $6.3715 \times 10^{-1}$  & $\phantom{-}1.2685\times 10^{-2}$  & $6.3725 \times 10^{-1}$   & $7$\\  
$4$   & $-1.0500\times 10^{-5}$   & $7.4620 \times 10^{-2}$  & $-5.3232\times 10^{-5}$    & $7.5483 \times 10^{-2}$    & $12$\\
$6$   & $-1.6859 \times 10^{-5}$  & $5.4567 \times 10^{-3}$  &  $-1.6849 \times 10^{-5}$  & $5.4592 \times 10^{-3}$   & $13$
\end{tabular}
\label{table:rail:loewner-vs-exact}
\end{table}
The decay of the singular values of $s_i\IL-\IM$ is shown in Figure \ref{fig:irka_loewner_decay}, indicating a natural choice for $k$. Only the $r=6$ case is presented; the other cases show the same pattern.

We compare the performance of \HinfIRKA~with that of \IRKA~for $r=2,4,6$, to illustrate the improvement offered  by optimizing over the $d_r$-term. The results are listed in Table \ref{rail_79841}. As expected,  \HinfIRKA~outperforms \IRKA~for every $r$-value. 
\begin{table}[htb]
\caption{Relative $\Hinf$ error  for the Cooling Steel Model, order $=79841$ }
\center
\begin{tabular}{ccccc}
$r$&\IRKA&\HinfIRKA\\
\hline
$2$&$6.46\times 10^{-1}$ &$\mathbf{6.37\times 10^{-1}}$ \\
$4$&$1.80 \times 10^{-1}$&$\mathbf{7.46 \times 10^{-2}}$\\
$6$ &$1.40 \times 10^{-2}$&$\mathbf{5.46\times 10^{-3}}$
\end{tabular}
\label{rail_79841}
\end{table}

\begin{figure}[hhht]
\centerline{\hbox{
\epsfysize=7cm
\epsfxsize=10cm
\epsffile{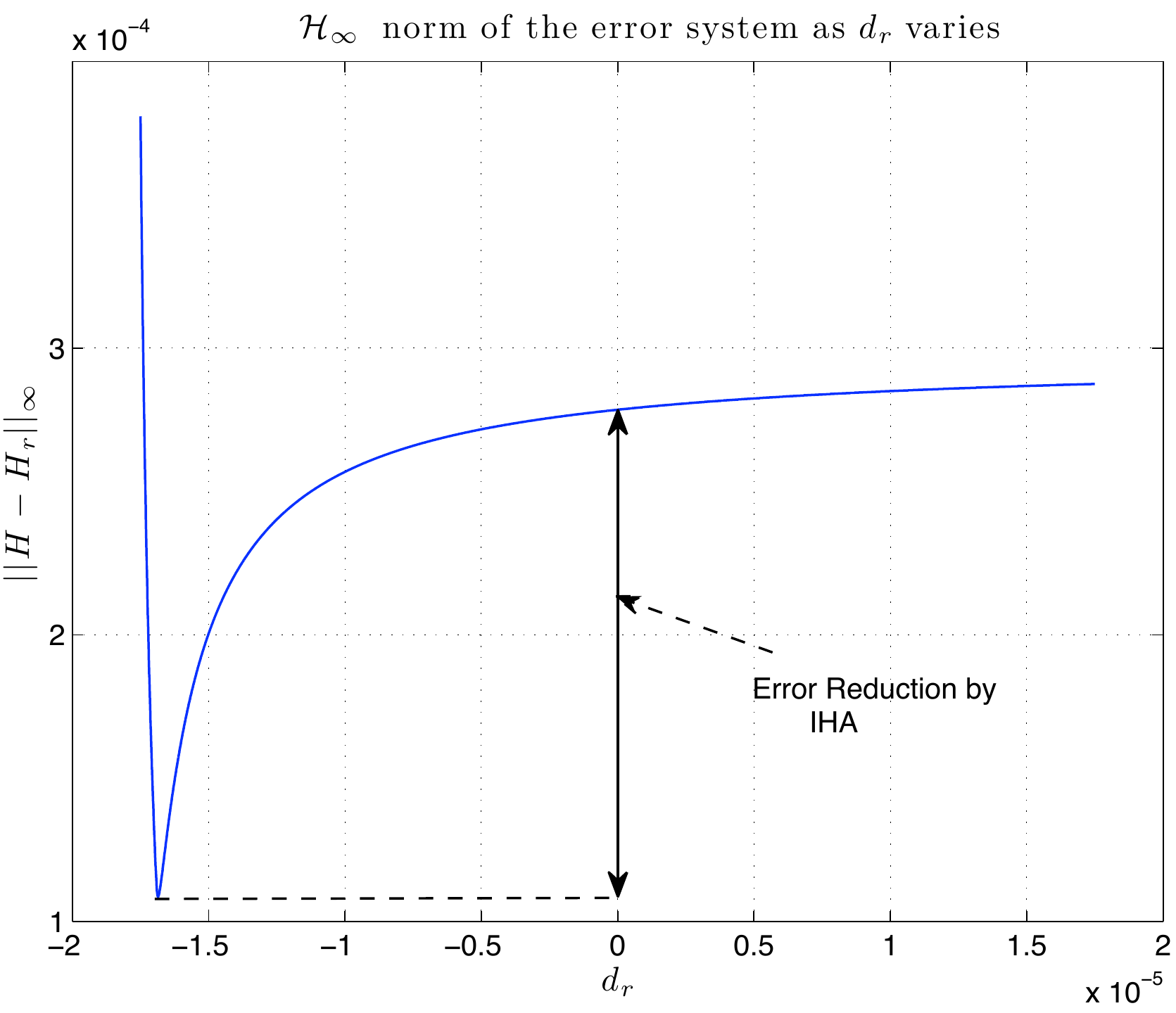}
}
}
\caption{ $\Hinf$ error as a function of the $\dr$-term}\label{rail_dr_gain}
\end{figure}

\begin{figure}[hhhht]
\centerline{\hbox{
\includegraphics[width=10.0cm]{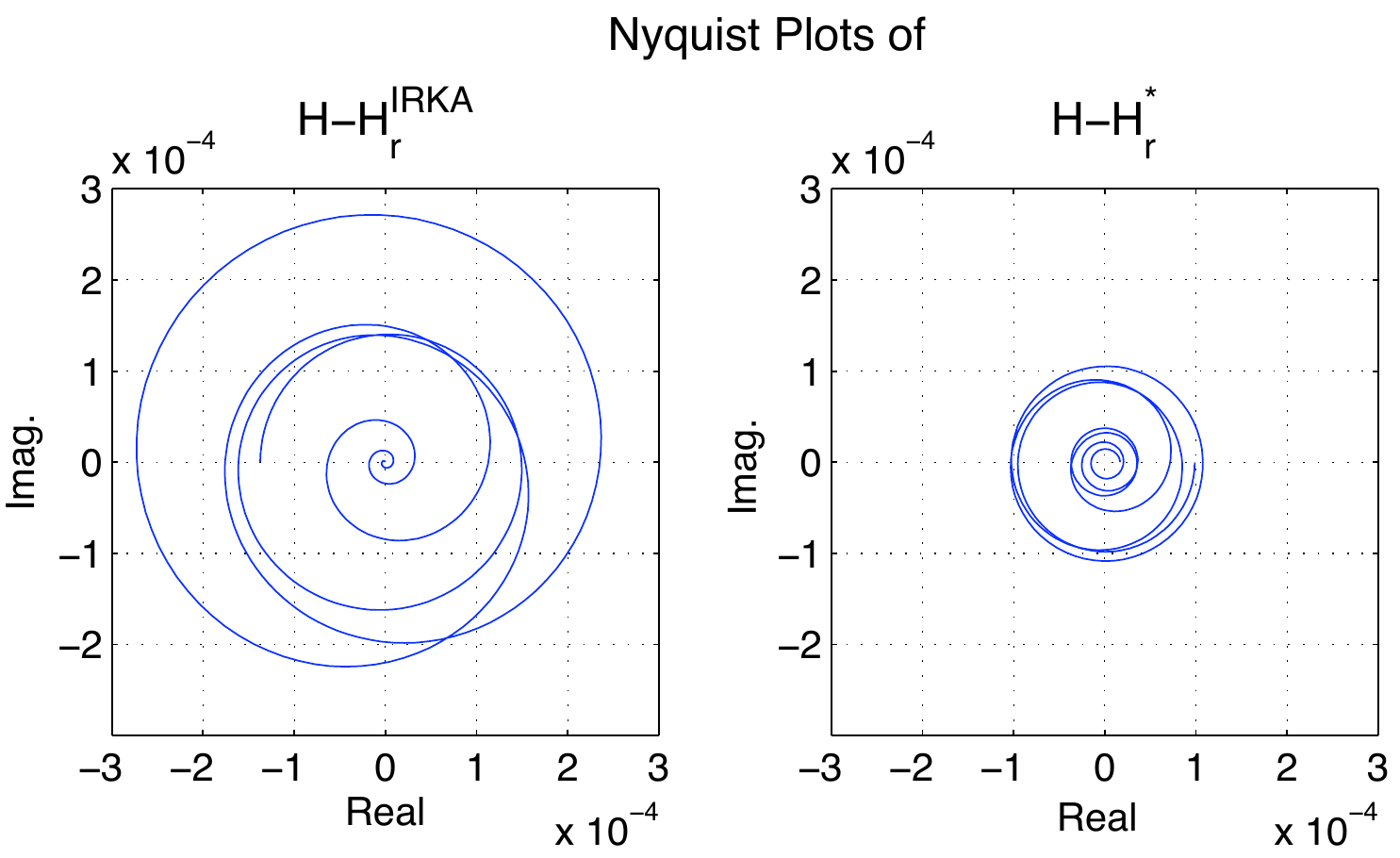}
}
}
\caption{ The error curve for $H_r^{\star}$ is nearly circular}\label{rail_nyquist}
\end{figure}

In Figure \ref{rail_dr_gain}, we 
demonstrate how the {\it absolute} $\Hinf$ error changes over values of $\dr$, for the order $r=6$ approximation.  In this case, the $\Hinf$-error is decreased by over a factor of two, from $2.58\times 10^{-4}$ down to $1.09 \times 10^{-4}$.  We also note that, for $r=6$, computing the optimal $\dr^{\star}$  term placed an additional interpolation point at $2.44\times10^{-2}$, yielding exactly $2r+1=13$ interpolation points in $\complex_{+}$ as suggested by Theorem \ref{thm:trefsufficient}.  Recalling Theorem \ref{thm:trefsufficient}, in Figure \ref{rail_nyquist} we demonstrate how this additional interpolation condition therefore results in a tighter, more circular Nyquist plot, which is equivalent to the image of the error along the imaginary axis being nearly circular.

\section{Conclusions}
We have introduced an interpolation-based model reduction technique to construct high-fidelity  $\Hinf$ approximations for large-scale linear dynamical systems.  For a given order $r$, $r$ Hermite interpolation points are produced that induce (locally) optimal $\Htwo$ approximation.
The $d_r$ (feed-forward) term is then adjusted in a such way that interpolation at  these initial $2r$ points is retained while an additional interpolation point is added that minimizes the $\Hinf$-error norm. By employing a data-driven Loewner approach, this last step may be performed at negligible cost; no large-scale $\Hinf$-norm computations are ever needed. 
The dominant cost of the method lies with the solution of sparse linear systems.
Four numerical examples show that the proposed method produces high fidelity $\Hinf$ reduced-order models that are better than those obtained by balanced truncation; and as good as (and sometimes better than) those obtained by optimal Hankel norm approximation; in all cases, at much lower computational cost. 

\bibliographystyle{plainnat}
\bibliography{references}

\end{document}